\documentclass[a4paper]{article}

\usepackage{graphicx}
\usepackage{amsmath}
\usepackage{amsthm}
\usepackage{amssymb, verbatim}
\usepackage[english]{babel}
\usepackage{mathrsfs}
\usepackage[compress]{cite}
\usepackage[usenames,dvipsnames]{color}
\usepackage{mdwlist}
\usepackage{paralist}

\renewcommand{\emph}{\textit}

\renewcommand{\underbar}{\underline}

\newtheorem{lemma}{\bf Lemma}[section]
\newtheorem{theorem}[lemma]{\bf Theorem}
\newtheorem{remark}[lemma]{\bf Remark}
\newtheorem{definition}[lemma]{Definition}
\newtheorem{assumption}[lemma]{Assumption}

\newtheorem{corollary}[lemma]{\bf Corollary}

\newcommand{\jpfig}[4]{\begin{figure}[t] \centering \includegraphics[width=#1\linewidth]{#2} \caption{\label{#3}#4} \end{figure}}
\renewcommand{\Re}{\mathbb{R}}
\renewcommand{\matrix}[2]{\left[\begin{array}{#1} #2 \end{array}\right] }
\newcommand{\diag}{\text{diag}}
\newcommand{\spanf}{\text{span}}
\newcommand{\nullf}{\text{null}}

\newcommand{\comm}{\mathcal{C}}
\newcommand{\p}{\mathcal{P}}
\newcommand{\A}{\mathcal{A}}
\newcommand{\B}{\mathcal{B}}
\newcommand{\K}{\mathcal{K}}

\newcommand{\IEEEQED}{~\rule[-1pt]{5pt}{5pt}\par\medskip}
\newenvironment{IEEEproof}{{\bf Proof:\ }}{ \hfill \IEEEQED}

\begin{document}

\title{Optimal Structured Static State-Feedback Control Design with Limited Model Information for Fully-Actuated Systems\thanks{ An early version of this paper is accepted for presentation at the American Control Conference, 2011 \cite{Farokhi10}.} \thanks{The work of F.~Farokhi and K.~H.~Johansson were supported by grants from the Swedish Research Council and the Knut and Alice Wallenberg Foundation. The work of C.~Langbort was supported, in part, by the US Air Force Office of Scientific Research (AFOSR) under grant number MURI FA 9550-10-1-0573.}}

\author{Farhad~Farokhi\thanks{F.~Farokhi and K.~H.~Johansson are with ACCESS Linnaeus Center, School of Electrical Engineering, KTH-Royal Institute of Technology, SE-100 44 Stockholm, Sweden. E-mails: \{farokhi,kallej\}@ee.kth.se },~C\'{e}dric~Langbort\thanks{C.~Langbort is with the Department of Aerospace Engineering and the Coordinated Science Laboratory, University of Illinois at Urbana-Champaign, Illinois, USA. E-mail: langbort@illinois.edu},~and~Karl~H.~Johansson$^\ddag$}

\date{}

\maketitle

\begin{abstract}
We introduce the family of limited model information control design methods, which construct controllers by accessing the plant's model in a constrained way, according to a given design graph. We investigate the closed-loop performance achievable by such control design methods for fully-actuated discrete-time linear time-invariant systems, under a separable quadratic cost.
We restrict our study to control design methods which produce structured static state feedback controllers, where each subcontroller can at least access the state measurements of
those subsystems that affect its corresponding subsystem. We compute the optimal control design strategy (in terms of the competitive ratio and domination metrics) when the control designer has access to the local model information and the global interconnection structure of the plant-to-be-controlled. Lastly, we study the trade-off between the amount of model information exploited by a control design method and the best closed-loop performance (in terms of the competitive ratio) of controllers it can produce.
\end{abstract}

\section{Introduction}
Many modern control systems, such as aircraft and satellite formation~\cite{Giulietti2000,Kapilal1999}, automated highways and other shared infrastructure~\cite{Swaroop1999,Negenborn2010}, flexible structures~\cite{Joshi1989}, and supply chains~\cite{Dunbar2007}, consist of a large number of subsystems coupled through their performance goals or system dynamics. When regulating this kind of plant, it is often advantageous to adopt a distributed control architecture, in which the controller itself is composed of interconnected subcontrollers, each of which accesses a strict subset of the plant's output. Several control synthesis methods have been proposed over the past decades that result in distributed controllers of this form, with various types of closed-loop stability and performance guarantees~(e.g.,~\cite{Levine1971,Wenk1980,Soderstrom1978,Shih-HoWang1973,Paganini1999,Bamieh2002,Lall2003,Hu1994,Scorletti2001}). Most recently, the tools presented in~\cite{Rotkowitz2006} and~\cite{Voulgaris200351} revealed how to exploit the specific interconnection of classes of plants (the so-called quadratically invariant systems) to formulate convex optimization problems for the design of structured $H_{\infty}$- and $H_2$- optimal controllers. A common thread in this part of the literature is the assumption that, even though the controller is structured, its design can be performed in a centralized fashion, with full knowledge of the plant model. However, in some applications (described in more detail in the next paragraph), this assumption is not always warranted, as the design of each subcontroller may need to be carried out by a different control designer, with no access to the global model of the plant, although its interconnection structure and the common closed-loop cost function to be minimized are public knowledge. This class of problems, which we refer to as ``limited model information control design problems'', is the main object of interest in the present paper.

Limited model information control design occurs naturally in contexts where the subsystems belong to different entities, which may consider their model information private and may thus be reluctant to share it with others. In this case, the designers may have to resort to ``communication-less'' strategies in which subcontroller $K_i$ depends solely on the description of subsystem $i$'s model. This case is well illustrated by supply chains, where the economic incentives of competing companies might limit the exchange of model information (such as, inventory volume, transportation efficiency, raw material sources, and decision process) inside a layer of the chain. Another reason for using communication-less strategies in more general design situations, even when the circulation of plant information is not restricted a priori, is that the resulting subcontroller $K_i$ does not need to be modified if the characteristics of a particular subsystem, which is not directly connected to subsystem~$i$, vary. For instance, consider a chemical plant in the process industry, with thousands of local controllers. In such a large-scale system, the tuning of each local controller should not require model parameters from other parts of the system so as to simplify maintenance and limit controller complexity. Note that engineers often implement these large-scale systems as a whole using commercially available pre-designed modules. These modules are designed, in advance, with no prior knowledge of their possible use or future operating condition. This lack of availability of the complete model of the plant, at the time of the design, constrains the designer to only use its own model parameters in each module's control design.

Control design based on uncertain plant model information is a classic topic in the robust control literature~\cite{Doyle1982,Zames1981,Ball1987,zhou1998}. However, designing an optimal controller without a global model is different from a robust control problem. In optimal control design with limited model information, subsystems do not have any prior information about the other subsystems' model; i.e., there is no nominal model for the design procedure and there is no bound on the model uncertainties. There have been some interesting approaches for tackling this problem. For instance, references~\cite{Ando1963,Sethi1998,Gajtsgori1979,Sezer1986} introduced methods for designing sub-optimal decentralized controllers without a global dynamical model of the system. In these papers, the authors assume that the large-scale system to be controlled consists of an interconnection of weakly coupled subsystems. They design an optimal controller for each subsystem using only the corresponding local model, and connect the obtained subcontrollers to construct a global controller. They show that, when coupling is negligible, this latter controller is satisfactory in terms of closed-loop stability and performance. However, as coupling strength increases, even closed-loop stability guarantees are lost. Other approaches such as~\cite{Dunbar2007,Negenborn2010} are based on receding horizon control and use decomposition methods to solve each step's optimization problem in a decentralized manner with only limited information exchange between subsystems. What is missing from the literature, however, is a rigorous characterization of the best closed-loop performance that can be attained through limited model information design and, a study of the trade off between the closed-loop performance and the amount of exchanged information. We tackle this question in the present paper for a particular class of systems (namely, the set of fully-actuated discrete-time linear time-invariant dynamical systems) and a particular class of control laws (namely, the set of structured linear static state feedback controllers where each subcontroller can at least access the state measurements of
those subsystems that affect its corresponding subsystem).

In this paper, we study the properties of limited model information control design methods. We investigate the relationship between the amount of plant information available to the designers, the nature of the plant interconnection graph, and the quality (measured by the closed-loop control goal) of controllers that can be constructed using their knowledge. To do so, we look at limited model information and communication-less control design methods as belonging to a special class of maps between the plant and controller sets, and make use of the competitive ratio and domination metrics introduced in~\cite{langbort10} to characterize their intrinsic limitations. To the best of our knowledge, there are no other metrics specifically tuned to control design methods. We address much more general classes of subsystems and of limitations on the model information available to the designer than is done in~\cite{langbort10}. Specifically, we consider limited model information structured static state-feedback control design for interconnections of fully-actuated (i.e., with invertible $B$-matrix) discrete-time linear time-invariant subsystems with quadratic separable (i.e., with block diagonal $Q$- and $R$-matrices) cost function. Our choice of such a cost function is motivated by our interest in applications such as power grids~\cite{Baughman1997,Berger1989,ChaoPeck1996,Botterud2005} and~\cite[Chs.~5,10]{Negenborn2010}, supply chains~\cite{Braun2003229,Dunbar2007}, and water level control~\cite[Ch.~18]{Negenborn2010}, which have been shown to be well-modeled by dynamically-coupled but cost-decoupled interconnected systems. We show in the last section of the paper that the assumption on the $B$-matrix can be partially removed for the sinks (i.e., subsystems that cannot affect any other subsystem) in the plant graph. 

We investigate the best closed-loop performance achievable by structured static state feedback controllers constructed by limited model information design strategies. We show that the result depends crucially on the plant graph and the control graph. In the case where the plant graph contains no sink and the control graph is a supergraph of the plant graph, we extend the fact proven in~\cite{langbort10} that the deadbeat strategy is the best communication-less control design method. However, the deadbeat control design strategy is dominated when the plant graph has sinks, and we exhibit a better, undominated, communication-less control design method, which, although having the same competitive ratio as the deadbeat control design strategy, takes advantage of the knowledge of the sinks' location to achieve a better closed-loop performance in average. We characterize the amount of model information needed to achieve better competitive ratio than the deadbeat control design strategy. This amount of information is expressed in terms of properties of the design graph; a directed graph which indicates the dependency of each subsystem's controller on different parts of the global dynamical model.

This paper is organized as follows. After formulating the problem of interest and defining the performance metrics in Section~\ref{sec:1}, we characterize the best communication-less control design method according to both competitive ratio and domination metrics in Section~\ref{sec:2}. In Section~\ref{sec:3}, we show that achieving a strictly better competitive ratio than these control design methods requires a complete design graph when the plant graph is itself complete. Finally, we end with a discussion on extensions in Section~\ref{subsec:6.1} and the conclusions in Section~\ref{sec:7}.

\subsection{Notation}
Sets will be denoted by calligraphic letters, such as $\mathcal{P}$ and $\mathcal{A}$. If $\mathcal{A}$ is a subset of $\mathcal{M}$ then $\mathcal{A}^c$ is the complement of $\mathcal{A}$ in $\mathcal{M}$, i.e., $\mathcal{M}\setminus \mathcal{A}$.

Matrices are denoted by capital roman letters such as $A$. $A_j$ will denote the $j^{\textrm{th}}$ row of $A$. $A_{ij}$ denotes a sub-matrix of matrix $A$, the dimension and the position of which will be defined in the text. The entry in the $i^{\textrm{th}}$ row and the $j^{\textrm{th}}$ column of the matrix $A$ is $a_{ij}$.

Let $S_{++}^n$ ($S_{+}^n$) be the set of symmetric positive definite (positive semidefinite) matrices in $\Re^{n\times n}$. $A > (\geq) 0$ means that the symmetric matrix $A\in \Re^{n\times n}$ is positive definite (positive semidefinite) and $A > (\geq) B$ means that $A-B > (\geq) 0$.

$\underline{\lambda}(Y)$ and $\bar{\lambda}(Y)$ denote the smallest and the largest eigenvalues of the matrix $Y$, respectively. Similarly, $\underline{\sigma}(Y)$ and $\bar{\sigma}(Y)$ denote the smallest and the largest singular values of the matrix $Y$, respectively. Vector $e_i$ denotes the column-vector with all entries zero except the $i^{\textrm{th}}$ entry, which is equal to one.

All graphs considered in this paper are directed, possibly with self-loops, with vertex set $\{1,...,q\}$ for some positive integer $q$. If $G=(\{1,...,q\},E)$ is a directed graph, we say that $i$ is a sink if there does not exist $j \neq i $ such that $(i,j) \in E$. A loop of length $t$ in $G$ is a set of distinct vertices $\{i_1,...,i_t\}$ such that $(i_t,i_1) \in E$ and $(i_p,i_{p+1}) \in E$ for all $ 1 \leq p \leq t-1$. We will sometimes refer to this loop as $(i_1 \rightarrow i_2 \rightarrow \dots \rightarrow i_t \rightarrow i_1)$. The adjacency matrix $S$ of graph $G$ is the $q \times q$ matrix whose entries satisfy
$$
s_{ij} = \left\{\begin{array}{cl} 1 & \mbox{ if } (j,i) \in E\\ 0 & \mbox{ otherwise.} \end{array} \right.
$$
Since the set of vertices is fixed here, a subgraph of $G$ is a graph whose edge set is a subset of the edge set of $G$ and a supergraph of $G$ is a graph of which $G$ is a subgraph. We use the notation $G'\supseteq G$ to indicate that $G'$ is a supergraph of $G$.

\section{Control Design with Limited Model Information} \label{sec:1}
 In this section, we introduce the system model and the problem under consideration, but first, we present a simple illustrative example.
\subsection{Illustrative Example}
Consider a discrete-time linear time-invariant dynamical system composed of three subsystems represented in state-space form as
\begin{equation*}
\begin{split}
\matrix{c}{x_1(k+1)\\x_2(k+1)\\x_3(k+1)}\hspace{-.05in}=\hspace{-.03in}&
\matrix{ccc}{a_{11}&a_{12}&0\\a_{21}&a_{22}&a_{23}\\0&a_{32}&a_{33}}
\matrix{c}{x_1(k)\\x_2(k)\\x_3(k)}\hspace{-.04in}+\hspace{-.04in}
\matrix{c}{b_{11}u_1(k)\\b_{22}u_2(k)\\b_{33}u_3(k)},
\end{split}
\end{equation*}
where, for each subsystem $i$, $x_i(k)\in\mathbb{R}$ is the state and $u_i(k)\in\mathbb{R}$ is the control signal. This system, which is illustrated in Figure~\ref{figure0_5}, is a simple networked control system. Networked control systems have several important characteristics. First, they are often distributed geographically. Therefore, it is natural to assume that a given subsystem can only influence its neighboring subsystems. We capture this fact using a directed graph called the plant graph like the one presented in Figure~\ref{figure0}($a$) for this example. This star graph corresponds to applications like unmanned aerial vehicles formation, platoon of vehicles, and composite formations of power systems~\cite{ChenStankovic2005,Zeynelgil2002}.

Second, any communication medium that we use to transmit the sensor measurements and actuation signals in networked control systems brings some limitations. For instance, every communication network has band-limited channels. Therefore, when designing subcontrollers, it might not make sense to assume that it can instantaneously access full state measurements of the plant. The state measurement availability in this example is
\begin{equation*}
\begin{split}
\matrix{c}{u_1(k)\\u_2(k)\\u_3(k)}=\matrix{ccc}{k_{11}&k_{12}&0\\k_{21}&k_{22}&k_{23}
\\0&k_{32}&k_{33}}\matrix{c}{x_1(k)\\x_2(k)\\x_3(k)}.
\end{split}
\end{equation*}
We use a control graph to characterize the controller structure. Control graph $G_\K$ in Figure~\ref{figure0}($b$) represents the state-measurement availability in this example. It corresponds to the case where neighboring subsystems transmit their state-measurements to each other, which is common for unmanned aerial vehicles formation, autonomous ground vehicles platoons, and biological system of particles~\cite{Giulietti2000,Kapilal1999,PhysRevLett751226,Jadbabaie1205192}.

Finally, in large-scale dynamical systems, it might be extremely difficult (if not impossible) to identify all system parameters and update them globally. One can only hope that the designer has access to the local parameter variations and update the corresponding subcontroller based on them. Therefore, it makes sense to assume that each local controller only has access to model information from its corresponding subsystem; i.e., designer of subcontroller $i$ uses only $\{a_{i1},a_{i2},a_{i3}\}$ in the design procedure
$$
[k_{i1}\;\;k_{i2}\;\;k_{i3}]=\Gamma_i\left([a_{i1}\;\;a_{i2}\;\;a_{i3}],b_{ii}\right),
$$
where $\Gamma_i:\mathbb{R}^3\times \mathbb{R}\rightarrow\mathbb{R}^3$ is the control design map. Note that the block-diagram in Figure~\ref{figure0_5} does not specify $\Gamma$. We will use a directed graph called the design graph to capture structural properties of $\Gamma$. In the rest of this section, we formalize the above notions for more general design problems.

\subsection{Plant Model}
Let a graph $G_{\mathcal{P}}=(\{1,...,q\},E_{\mathcal{P}})$ be given, with adjacency matrix $S_{\mathcal{P}} \in \{0,1\}^{q \times q}$. We define the following set of matrices associated with $S_{\mathcal{P}}$:
\begin{equation}
\begin{split}
\mathcal{A}(S_{\mathcal{P}})=\{ A \in \; \mathbb{R}^{n \times n} \;|\; A_{ij} = 0& \in \mathbb{R}^{n_i \times n_j} \mbox{ for all } \\ & 1\leq i,j \leq q \mbox { such that } (s_{\mathcal{P}})_{ij}=0 \},
\end{split}
\end{equation}
where for each $1 \leq i \leq q$, integer number $n_i$ is the dimension of subsystem $i$. Implicit in these definitions is the fact that $\sum_{i=1}^q n_i=n$. Also, for a given scalar $\epsilon >0$, we let
\begin{equation} \label{eqn:defBe}
\begin{split}
\mathcal{B}(\epsilon) =\{ B \in \mathbb{R}^{n \times n} \; | &\; \underbar{\sigma}(B) \geq \epsilon, B_{ij} = 0 \in \mathbb{R}^{n_i \times n_j} \mbox{ for all } 1 \leq i\neq j \leq q\}.
\end{split}
\end{equation}
The set $\mathcal{B}(\epsilon)$ defined in~(\ref{eqn:defBe}) is made of invertible block-diagonal square matrices since $\underbar{\sigma}(B) \geq \epsilon>0$ for each matrix $B\in\B(\epsilon)\subseteq\mathbb{R}^{n\times n}$. With these definitions, we can introduce the set $\mathcal{P}$ of plants of interest as the space of all discrete-time linear time-invariant dynamical systems of the form
\begin{equation}
x(k+1) = Ax(k) + Bu(k) \; ; \; x(0)=x_0,
\label{firsteq}
\end{equation}
with $A \in \mathcal{A}(S_{\mathcal{P}})$, $B \in \mathcal{B}(\epsilon)$, and $x_0 \in \mathbb{R}^n$. Clearly $\mathcal{P}$ is isomorph to $\mathcal{A}(S_{\mathcal{P}}) \times \mathcal{B}(\epsilon) \times \mathbb{R}^n$ and, slightly abusing notation, we will thus identify a plant $P \in \mathcal{P}$ with the corresponding triple $(A,B,x_0)$.

A plant $P \in \mathcal{P}$ can be thought of as the interconnection of $q$ subsystems, with the structure of the interconnection specified by the graph $G_{\mathcal{P}}$ (i.e., subsystem $j$'s output feeds into subsystem $i$ only if $(j,i) \in E_{\mathcal{P}}$). As a consequence, we refer to $G_{\mathcal{P}}$ as the ``plant graph''. We will denote the ordered set of state indices pertaining to subsystem $i$ as $\mathcal{I}_i$; i.e., $\mathcal{I}_i:=(1+\sum_{j=1}^{i-1} n_j, \dots, n_i+\sum_{j=1}^{i-1} n_j)$. For subsystem $i$, state vector and input vector are defined as
\begin{equation*}
\underbar{x}_i = \left[x_{\ell_1}\;\cdots\;x_{\ell_{n_i}}\right]^T, \hspace{.2in} \underbar{u}_i = \left[u_{\ell_1}\;\cdots\;u_{\ell_{n_i}} \right]^T
\end{equation*}
where the ordered set of indices $(\ell_1,\dots,\ell_{n_i})\equiv\mathcal{I}_i$, and its dynamics is specified by
$$
\underbar{x}_i(k+1) = \sum_{j=1}^q A_{ij} \underbar{x}_j(k) + B_{ii} \underbar{u}_i(k).
$$
According to the specific structure of $\B(\epsilon)$ given in~(\ref{eqn:defBe}), each subsystem is fully-actuated, with as many input as states, and controllable in one time-step. Possible generalization of the results to a (restricted) family of under-actuated systems is discussed in Section~\ref{subsec:6.1}.

\begin{figure}[t] \centering \includegraphics[width=0.7\linewidth]{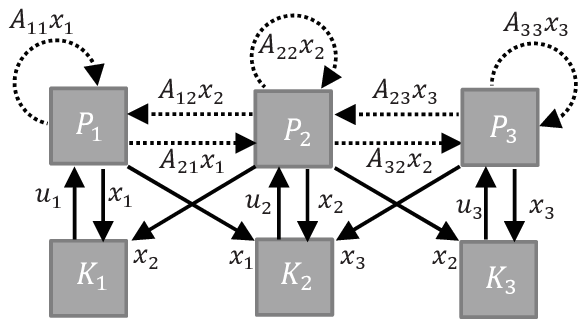} \caption{\label{figure0_5} Physical interconnection between different subsystems and controllers corresponding to $G_\p$ and $G_\K$ in Figures~\ref{figure0}($a$)~and~\ref{figure0}($b$), respectively.} \vspace{.1in} \includegraphics[width=0.5\linewidth]{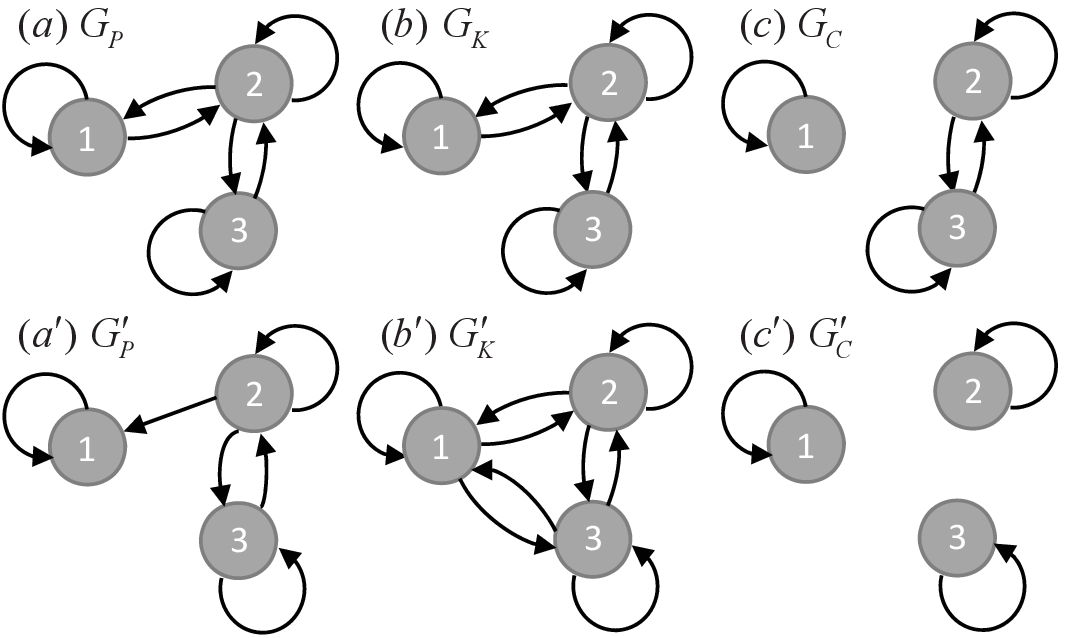} \caption{\label{figure0} $G_\p$ and $G'_\p$ are examples of plant graphs, $G_\K$ and $G'_\K$ are examples of control graphs, and $G_\comm$ and $G'_\comm$ are examples of design graphs.}
\end{figure}

Figure~\ref{figure0}($a$) shows an example of a plant graph $G_p$. Each node represents a subsystem of the system. For instance, the second subsystem in this example may affect the first subsystem and the third subsystem; i.e., sub-matrices $A_{12}$ and $A_{32}$ can be nonzero. The self-loop for the second subsystem shows that $A_{22}$ may be non-zero. The plant graph $G_\p$ in Figure~\ref{figure0}($a$) does not contain any sink. In contrast,
the first subsystem of the plant graph $G'_\p$ in Figure~\ref{figure0}($a'$) is a sink. The control graph $G_\K$ is introduced in the next subsection.

\subsection{Controller Model}
Let a control graph $G_{\mathcal{K}}$ be given, with adjacency matrix $S_{\mathcal{K}}$. The control laws of interest in this paper are linear static state-feedback control laws of the form
\begin{equation*}
u(k)=Kx(k),
\end{equation*}
where
\begin{equation}
\begin{split}
K\in \mathcal{K}(S_{\mathcal{K}})=\{ K \in \; \mathbb{R}^{n \times n} | K_{ij} = 0& \in \mathbb{R}^{n_i \times n_j} \mbox{ for} \\ \mbox{all }& 1\leq i,j \leq q \mbox { such that } (s_{\mathcal{K}})_{ij}=0 \}.
\end{split}
\end{equation}
In particular, when $G_{\mathcal{K}}$ is a complete graph, $\mathcal{K}(S_{\mathcal{K}}) = \mathbb{R}^{n \times n}$, while, if $G_{\mathcal{K}}$ is totally disconnected with self-loops, $\mathcal{K}(S_{\mathcal{K}})$ represents the set of fully-decentralized controllers. When adjacency matrix $S_{\mathcal{K}}$ is not relevant or can be deduced from context, we refer to the set of controllers as $\mathcal{K}$.

An example of a control graph $G_\K$ is given in Figure~\ref{figure0}($b$). Each node represents a subsystem-controller pair of the overall system. For instance, Figure~\ref{figure0}($b$) shows that the first subsystem's controller can use state measurements of the second subsystem besides its own state measurements. Figure~\ref{figure0}($b'$) shows a complete graph, which indicates that each subsystem has access to full state measurements of all other subsystems; i.e., $\mathcal{K}(S_{\mathcal{K}})=\mathbb{R}^{n \times n}$.

\subsection{Linear State Feedback Control Design Methods}
A control design method $\Gamma$ is a map from the set of plants $\mathcal{P}$ to the set of controllers $\mathcal{K}$. Just like plants and controllers, a control design method can exhibit structure which, in turn, can be captured by a design graph. Let a control design method $\Gamma$ be partitioned according to subsystems dimensions as
\begin{equation}
\Gamma=\matrix{ccc}{ \Gamma_{11} & \cdots & \Gamma_{1q} \\ \vdots & \ddots & \vdots \\ \Gamma_{q1} & \cdots & \Gamma_{qq} }
\label{eq_gamma}
\end{equation}
and a graph $G_{\mathcal{C}}=(\{1,...,q\}, E_{\mathcal{C}})$ be given, with adjacency matrix $S_{\mathcal{C}}$. Each block $\Gamma_{ij}$ represents a map $\mathcal{A}(S_{\mathcal{P}}) \times \mathcal{B}(\epsilon) \rightarrow \mathbb{R}^{n_i \times n_j}$. Control design method $\Gamma$ can be further partitioned in the form
$$ \Gamma=\matrix{ccc}{ \gamma_{11} & \cdots & \gamma_{1n} \\ \vdots & \ddots & \vdots \\ \gamma_{n1} & \cdots & \gamma_{nn} }, $$
where each $\gamma_{ij}$ is a map $\mathcal{A}(S_{\mathcal{P}}) \times \mathcal{B}(\epsilon) \rightarrow \mathbb{R}$. We say that $\Gamma$ has structure $G_{\mathcal{C}}$ if, for all $i$, the map $\left[\Gamma_{i1} \; \cdots \; \Gamma_{iq} \right]$ is only a function of
\begin{equation}
\left\{ \left[A_{j1}\;\cdots\;A_{jq}\right], B_{jj} \; | \; (s_{\mathcal{C}})_{ij} \neq 0 \right \}.
\end{equation}
In words, a control design method has structure $G_{\mathcal{C}}$ if and only if, for all $i$, the subcontroller of subsystem $i$ is constructed with knowledge of the plant model of only those subsystems $j$ such that $(j,i) \in E_{\mathcal{C}}$. The set of all control design methods with structure $G_{\mathcal{C}}$ will be denoted by $\mathcal{C}$. In the particular case where $G_{\mathcal{C}}$ is the totally disconnected graph with self-loops (meaning that every node in the graph has a self-loop; i.e, $S_{\mathcal{C}} = I_q$), we say that a control design method in $\mathcal{C}$ is ``communication-less'', so as to capture the fact that subsystem $i$'s subcontroller is constructed with no information coming from (and, hence, no communication with) any other subsystem $j$, $j \neq i$. Therefore, the design graph indicates knowledge (or lack thereof) of entire block rows in the aggregate system matrix. When $G_{\mathcal{C}}$ is not a complete graph, we refer to $\Gamma \in \mathcal{C}$ as being ``a limited model information control design method''.

Note that $\mathcal{C}$ can be considered as a subset of the set of functions from $\mathcal{A}(S_{\mathcal{P}}) \times \mathcal{B}(\epsilon)$ to $\mathcal{K}(S_{\mathcal{K}})$, since a design method with structure $G_{\mathcal{C}}$ is not a function of initial state $x_0$. Hence, when $\Gamma \in \mathcal{C}$ we will write $\Gamma(A,B)$ instead of $\Gamma(P)$ for plant $P=(A,B,x_0) \in \mathcal{P}$.

An example of a design graph $G_{\mathcal{C}}$ is given in Figure~\ref{figure0}($c$). Each node represents a subsystem-controller pair of the overall system. For instance, $G_\comm$ shows that the third subsystem's model is available to the designer of the second subsystem's controller but not the first subsystem's model. Figure~\ref{figure0}($c'$) shows a fully disconnected design graph with self-loops $G'_\comm$. A local designer in this case can only rely on the model of its corresponding subsystem; i.e., the design strategy is communication-less.

\subsection{Performance Metrics} \label{subsec:PM}
The goal of this paper is to investigate the influence of the plant and design graph on the properties of controllers constructed by limited model information control design methods. To this end, we will use two performance metrics for control design methods. These performance metrics are adapted from the notions of competitive ratio and domination introduced in~\cite{langbort10}, so as to take plant, controller, and control design structures into account. Following the approach in~\cite{langbort10}, we start by associating a closed-loop performance criterion to each plant $P=(A,B,x_0) \in \p$ and controller $K\in \K$. As explained in the introduction, we are particularly interested in dynamically-coupled but cost-decoupled systems in this paper, hence, we use a cost of the form
\begin{equation} \label{eqn:1}
J_P (K)=\sum_{k=1}^\infty x(k)^TQx(k) +\sum_{k=0}^\infty u(k)^TRu(k),
\end{equation}
where $Q \in S_{++}^n$ and $R \in S_{++}^n$ are block diagonal matrices, with each diagonal block entry belonging to~$S_{++}^{n_i}$. Note that the summation in the first term on the right-hand side of~(\ref{eqn:1}) starts from $k=1$. This is without loss of generality as the removed term $x(0)^TQx(0)$ is not a function of the controller. We make the following two standing assumptions:

\begin{assumption} $Q=R=I$. \end{assumption}

This is without loss of generality because the change of variables $(\bar{x},\bar{u})= (Q^{1/2}x,R^{1/2}u)$ transforms the performance criterion and state space representation
into
\begin{equation}
\label{cost_easy}
J_P (K)=\sum_{k=1}^\infty \bar{x}(k)^T\bar{x}(k)+\sum_{k=0}^\infty \bar{u}(k)^T\bar{u}(k),
\end{equation}
and
\begin{equation*}
\begin{split}
\bar{x}(k+1)&=Q^{1/2}AQ^{-1/2}\bar{x}(k)+Q^{1/2}BR^{-1/2}\bar{u}(k)
\\&=\bar{A}\bar{x}(k)+\bar{B}\bar{u}(k),
\end{split}
\end{equation*}
respectively, without affecting the plant, control, or design graph (due to the block diagonal structure of $Q$ and~$R$).

\begin{assumption} The set of matrices $\mathcal{B}(\epsilon)$ is replaced with the set of diagonal matrices with diagonal entries greater than or equal to $\epsilon$. \end{assumption}

This assumption is without loss of generality. Indeed, consider a plant $P=(A,B,x_0)\in \p$. Every sub-system's $B_{ii}$ matrix has a singular value decomposition $B_{ii}=U_{ii} \Sigma_{ii} V_{ii}^T$ with $\Sigma_{ii}\geq \epsilon I_{n_i\times n_i}$. Combining these singular value decompositions together results in a singular value decomposition for matrix $B=U\Sigma V^T$ where $U=\diag (U_{11},U_{22},\cdots,U_{qq})$, $\Sigma=\diag (\Sigma_{11},\Sigma_{22},\cdots,\Sigma_{qq})$, and $V=\diag (V_{11},V_{22},\cdots,V_{qq})$. Defining $\bar{x}(k)=U^Tx(k)$ and $\bar{u}(k)=V^Tu(k)$ results in
\begin{equation*}
\bar{x}(k+1)=U^TAU\bar{x}(k)+U^TBV\bar{u}(k),
\end{equation*}
where $U^TBV$ is diagonal. Because of the block diagonal structure of matrices $U$ and $V$, the change of variables $(A,B,x_0) \mapsto (U^TAU,U^TBV,U^Tx_0)$ does not affect the plant, control, or design graph. In addition, the cost function becomes
\begin{equation*}
\begin{split}
J_P (K)&=\sum_{k=1}^\infty \bar{x}(k)^TU^TU\bar{x}(k) +\sum_{k=0}^\infty \bar{u}(k)^TV^TV\bar{u}(k)\\&=\sum_{k=1}^\infty \bar{x}(k)^T\bar{x}(k)
 +\sum_{k=0}^\infty \bar{u}(k)^T\bar{u}(k),
\end{split}
\end{equation*}
which is of the form~(\ref{cost_easy}), because both $U$ and $V$ are unitary matrices. We are now ready to define the performance metrics of interest in this paper.

\begin{definition}\textit{(Competitive Ratio)} Let a plant graph $G_{\mathcal{P}}$, control graph $G_{\mathcal{K}}$ and constant $\epsilon > 0$ be given. Assume that, for every plant $P \in \mathcal{P}$, there exists an optimal controller $K^*(P) \in \mathcal{K}$ such that
\begin{equation*}
J_P (K^*(P))\leq J_P (K), \hspace{0.05in} \forall K \in \K.
\end{equation*}
The competitive ratio of a control design method $\Gamma$ is defined as
\begin{equation*}
r_{\p} (\Gamma)= \sup_{P=(A,B,x_0) \in \p} \frac{J_P (\Gamma(A,B))}{J_P (K^*(P))},
\end{equation*}
with the convention that ``$\frac{0}{0}$'' equals one.
\label{def_comp_rat}
\end{definition}
Note that the mapping $ K^* : P \rightarrow K^*(P)$ is not itself required to lie in the set $\comm$, as every component of the optimal controller may depend on all entries of the model matrices $A$ and $B$.

\begin{definition}\textit{(Domination)} A control design method $\Gamma$ is said to dominate another control design method $\Gamma'$ if
\begin{equation}
J_P(\Gamma(A,B))\leq J_P(\Gamma'(A,B)),\hspace{0.1in} \forall \; P=(A,B,x_0)\in \p,
\label{comp}
\end{equation}
with strict inequality holding for at least one plant in $\p$. When $\Gamma' \in \mathcal{C}$ and no control design method $\Gamma \in \mathcal{C}$ exists that satisfies~(\ref{comp}), we say that $\Gamma'$ is undominated in $\mathcal{C}$ for plants in $\mathcal{P}$.
\end{definition}

\subsection{Problem Formulation}
With the definitions of the previous subsections in hand, we can reformulate the main question of this paper regarding the connection between closed-loop performance, plant structure, and limited model information control design as follows. For a given plant graph, control graph, and design graph, we would like to determine
\begin{equation} \label{eqn:0}
\arg \min_{\Gamma \in \comm} r_{\p} (\Gamma).
\end{equation}
Since several design methods may achieve this minimum, we are interested in determining which ones of these strategies are \textit{undominated}.

In~\cite{langbort10}, this problem was solved in the case when $G_{\mathcal{P}}$ and $G_{\mathcal{K}}$ are complete graphs, $G_{\mathcal{C}}$ is a totally disconnected graph with self-loops (i.e., $S_{\mathcal{C}}=I_{q}$), and $\mathcal{B}(\epsilon)$ is replaced with singleton $\{I_n\}$. In this paper, we investigate the role of more general plant and design graphs. We also extend the results in~\cite{langbort10} for scalar subsystems to subsystems of arbitrary order $n_i\geq 1$, $1\leq i\leq q$.

\section{Plant Graph Influence on Achievable Performance} \label{sec:2}
\setcounter{subsection}{0}
In this section, we study the relationship between the plant graph and the achievable closed-loop performance in terms of the competitive ratio and domination.
\begin{definition} The deadbeat control design method $\Gamma^{\Delta}: \A(S_{\mathcal{P}}) \times \B(\epsilon) \rightarrow \K$ is defined as
\begin{equation*}
\Gamma^{\Delta} (A,B)=-B^{-1}A, \; \textrm{for all} \hspace{0.05in} P=(A,B,x_0)\in \p.
\end{equation*}
\end{definition}
This control design method is communication-less; i.e., the control design for the subsystem $i$ is a function of the model of subsystem $i$ only, because subsystem $i$'s controller gain $\left[ \Gamma^\Delta_{i1} (A,B) \; \cdots \; \Gamma^\Delta_{iq} (A,B) \right]$ equals to $B_{ii}^{-1} \left[ A_{i1} \; \cdots \; A_{iq} \right]$. The name ``deadbeat'' comes from the fact that the closed-loop system obtained by applying controller $\Gamma^{\Delta}(A,B)$ to plant $P=(A,B,x_0)$ reaches the origin in just one time-step~\cite{Emami-Naeini1982}.

\begin{remark} Note that for the case where the control graph $G_{\mathcal{K}}$ is a complete graph; i.e., $\mathcal{K}=\mathbb{R}^{n \times n}$, there exists a controller $K^*(P)$ satisfying the assumptions of Definition~\ref{def_comp_rat} for all $P \in \mathcal{P}$, namely, the optimal linear quadratic regulator which is independent of the initial condition of the plant. For incomplete control graphs, the optimal control design strategy $K^*(P)$ (if exists) might become a function of the initial condition~\cite{Levine1970}. Hence, we will use $K^*(A,B)$ instead of $K^*(P)$ when the control graph $G_{\mathcal{K}}$ is a complete graph for each plant $P=(A,B,x_0)\in\p$ to emphasize this fact. \end{remark}

 From Definition~\ref{def_comp_rat}, the notation $K^*(P)$ is reserved for the optimal control design strategy for any given control graph $G_\K$. In contrast, when $G_\K$ is not the complete graph, we will refer to the optimal \textit{unstructured} controller as $K^*_C(A,B)$.

\begin{lemma} \label{lem:0} Let the control graph $G_{\mathcal{K}}$ be a complete graph. The cost of the optimal control design strategy $K^*$ is lower-bounded by
\begin{equation*}
J_P(K^*(A,B)) \geq \left( \frac{\underline{\sigma}^2(B)}{\underline{\sigma}^2(B)+1} \right) J_P(\Gamma^\Delta(A,B)),
\end{equation*}
for all plants $P=(A,B,x_0)\in \p$.
\end{lemma}
\begin{IEEEproof} See Appendix~\ref{prooflem:0}. \end{IEEEproof}

\begin{theorem} \label{tho:1} Let the plant graph $G_{\mathcal{P}}$ contain no isolated node and $G_\K\supseteq G_\p$. Then the competitive ratio of the deadbeat control design method $\Gamma^\Delta$ is
$$
r_{\mathcal{P}}(\Gamma^{\Delta})=1+1/\epsilon^2.
$$
\end{theorem}
\begin{IEEEproof} Irrespective of the control graph $G_\K$ and for all plants $P\in\p$, it is true that $J_P(K^*_C(A,B)) \leq J_P(K^*(P))$. Therefore, we get
\begin{equation} \label{eqn:KC*K*}
\begin{split}
\frac{J_P(\Gamma^\Delta(A,B))}{J_P(K^*(P))}\leq \frac{J_P(\Gamma^\Delta(A,B))}{J_P(K^*_C(A,B))}.
\end{split}
\end{equation}
Now, using Lemma~\ref{lem:0}, we know that
\begin{equation} \label{eqn:DeltaK*C}
\begin{split}
\frac{J_P(\Gamma^\Delta(A,B))}{J_P(K^*_C(A,B))} \leq 1+\frac{1}{\underline{\sigma}^2(B)},
\end{split}
\end{equation}
for all $P=(A,B,x_0) \in \mathcal{P}$. Combining~(\ref{eqn:DeltaK*C}) and~(\ref{eqn:KC*K*}) results in
\begin{equation*}
\begin{split}
r_\p(\Gamma^\Delta)&=\sup_{P\in \p} \frac{J_P(\Gamma^\Delta(A,B))}{J_P(K^*(P))} \leq 1+\frac{1}{\epsilon^2}.
\end{split}
\end{equation*}
To show that this upper bound is attained, let us pick $i_1\in \mathcal{I}_{i}$ and $j_1\in \mathcal{I}_{j}$ where $1\leq i\neq j \leq q$ and $(s_{\mathcal{P}})_{ij}\neq 0$ (such indices $i$ and $j$ exist because plant graph $G_{\mathcal{P}}$ has no isolated node by assumption). Consider the system $A=e_{i_1}e_{j_1}^T$ and $B=\epsilon I$. The unique positive definite solution of the discrete algebraic Riccati equation
\begin{equation} \label{eqn:Riccati}
A^TXA-A^TXB(I+B^TXB)^{-1}B^TXA=X-I,
\end{equation}
is $X=I+[1/(1+\epsilon^2)]e_{j_1}e_{j_1}^T$. Consequently, the centralized controller $K^*_C(A,B)=-\epsilon/(1+\epsilon^2)e_{i_1}e_{j_1}^T$ belongs to the set $\K(S_\K)$ because $G_\K\supseteq G_\p$. Thus, we get
\begin{equation} \label{eqn:ineq11}
J_{(A,B,e_{j_1})}(K^*(A,B,e_{j_1})) \leq J_{(A,B,e_{j_1})}(K^*_C(A,B))
\end{equation}
since $K^*(P)$ has a lower cost than any other controller in $\K(S_\K)$. On the other hand, it is evident that
\begin{equation} \label{eqn:ineq21}
J_{(A,B,e_{j_1})}(K^*_C(A,B)) \leq J_{(A,B,e_{j_1})}(K^*(A,B,e_{j_1}))
\end{equation}
because the centralized controller has access to more state measurements. Using~(\ref{eqn:ineq11}) and~(\ref{eqn:ineq21}) simultaneously results in
\begin{equation*}
\begin{split}
J_{(A,B,e_{j_1})}(K^*(A,B,e_{j_1}))&=J_{(A,B,e_{j_1})}(K^*_C(A,B))\\&=1/(1+\epsilon^2).
\end{split}
\end{equation*}
 On the other hand $\Gamma^\Delta(A,B)=-[1/\epsilon] e_{i_1}e_{j_1}^T$ and $J_{(A,B,e_{j_1})}(\Gamma^\Delta(A,B))=1/\epsilon^2$. Therefore, $r_{\p}(\Gamma^{\Delta})=1+1/\epsilon^2$.
\end{IEEEproof}

\begin{remark} Consider the limited model information design problem given by the plant graph $G_\p$ in Figure~\ref{figure0}($a$) and the control graph $G'_\K$ in Figure~\ref{figure0}($b'$). Theorem~\ref{tho:1} shows that, if we apply the deadbeat control design strategy to this particular problem, the performance of the deadbeat control design strategy, at most, can be $1+1/\epsilon^2$ times the cost of the optimal control design strategy $K^*$. For instance, when $\B=\{I\}$ as in~\cite{langbort10}, we have $1+1/\epsilon^2=2$ since in this case $\epsilon=1$. Therefore, the deadbeat control design strategy is never worse than twice the optimal controller in this case.
\end{remark}

\begin{remark} There is no loss of generality in assuming that there is no isolated node in the plant graph $G_{\mathcal{P}}$, since it is always possible to design a controller for an isolated subsystem without any model information about the other subsystems and without impacting cost~(\ref{eqn:1}). In particular, this implies that there are $q \geq 2$ vertices in the graph because for $q=1$ the only subsystem that exists is an isolated node in the plant graph.
\end{remark}

\begin{remark} For implementation of the deadbeat control design strategy in each node, we only need the state measurements of the neighbors of that node. For the implementation of the optimal control design strategy $K^*$ when the control graph has many more links than the plant graph, the controller gain $K^*(P)$ is not necessarily a sparse matrix.
\end{remark}

With this characterization of $\Gamma^{\Delta}$ in hand, we are now ready to tackle problem~(\ref{eqn:0}).

\subsection{First case: plant graph $G_{\mathcal{P}}$ with no sink}
In this subsection, we show that the deadbeat control method $\Gamma^\Delta$ is undominated by communication-less control design methods for plants in $\mathcal{P}$, when $G_\mathcal{P}$ contains no sink. We also show that $\Gamma^\Delta$ exhibits the smallest possible competitive ratio among such control design methods. First, we state the following two lemmas.
\begin{lemma} \label{lem:1} Let the plant graph $G_{\mathcal{P}}$ contain no isolated node, the design graph $G_{\mathcal{C}}$ be a totally disconnected graph with self-loops, and $G_\K\supseteq G_\p$. A control design method $\Gamma \in \comm$ has bounded competitive ratio only if the following implication holds for all $1\leq i \leq q$ and all $j$:
$$
a_{\ell j}=0 \mbox{ for all } \ell \in \mathcal{I}_i \Rightarrow \gamma_{\ell j}(A,B)=0 \mbox{ for all } \ell \in \mathcal{I}_i,
$$
where $\mathcal{I}_i$ is the set of indices related to subsystem $i$; i.e., $\mathcal{I}_i=(1+\sum_{z=1}^{i-1} n_z, \dots, n_i+\sum_{z=1}^{i-1} n_z)$.
\end{lemma}
\begin{IEEEproof} See Appendix~\ref{prooflem:1}.\end{IEEEproof}

\begin{lemma} \label{lem:2} Let the plant graph $G_{\mathcal{P}}$ contain no isolated node, the design graph $G_{\mathcal{C}}$ be a totally disconnected graph with self-loops, and $G_\K\supseteq G_\p$. Assume the plant graph $G_\p$ has at least one loop. Then,
\begin{equation}
\label{eq_rat}
r_{\p}(\Gamma) \geq 1+ 1/\epsilon^2
\end{equation}
for all limited model information control design method $\Gamma$ in $\comm$. \end{lemma}
\begin{IEEEproof} See Appendix~\ref{prooflem:2}.\end{IEEEproof}

Using these two lemmas, we are ready to state and prove one of the main theorems in this paper and, as a result, find the solution to problem~(\ref{eqn:0}) when the plant graph $G_\mathcal{P}$ contains no sink.

\begin{theorem} \label{tho:3} Let the plant graph $G_{\mathcal{P}}$ contain no isolated node and no sink, the design graph $G_{\mathcal{C}}$ be a totally disconnected graph with self-loops, and $G_\K \supseteq G_\p$. Then the competitive ratio of any control design strategy $\Gamma\in\comm$ satisfies
$$
r_{\p} (\Gamma)\geq 1+1/\epsilon^2.
$$
\end{theorem}
\begin{IEEEproof} From Lemma~1.4.23 in~\cite{graphtheory}, we know that a directed graph with no sink must have at least one loop. Hence $G_{\mathcal{P}}$ must contain a loop. The result then follows from Lemma~\ref{lem:2}.\end{IEEEproof}

\begin{remark} Theorem~\ref{tho:3} shows that $r_{\p} (\Gamma)\geq r_{\p} (\Gamma^\Delta)$ for any control design strategy $\Gamma\in\comm$, and as a result the deadbeat control design method $\Gamma^{\Delta}$ becomes a minimizer of the competitive ratio function $r_{\mathcal{P}}$ over the set of communication-less design methods.
\end{remark}

We now turn our attention to domination properties of the deadbeat control design strategy.

\begin{lemma} \label{lem:undomination} Let the plant graph $G_{\mathcal{P}}$ contain no isolated node, the design graph $G_{\mathcal{C}}$ be a totally disconnected graph with self-loops, and $G_\K\supseteq G_\p$. The deadbeat control design strategy $\Gamma^\Delta$ is undominated, if there is no sink in the plant graph $G_\p$.
\end{lemma}
\begin{IEEEproof} See Appendix~\ref{prooflem:undomination}. \end{IEEEproof}

The following theorem shows that the deadbeat control design strategy is undominated by communication-less design methods if and only if the plant graph $G_{\mathcal{P}}$ has no sink. It thus provides a good trade-off between worst-case and average performance.

\begin{theorem} \label{tho:2} Let the plant graph $G_{\mathcal{P}}$ contain no isolated node, the design graph $G_{\mathcal{C}}$ be a totally disconnected graph with self-loops, and $G_\K\supseteq G_\p$. Then the deadbeat control design method $\Gamma^\Delta$ is undominated in $\mathcal{C}$ for plants in $\mathcal{P}$ if and only if the plant graph $G_{\mathcal{P}}$ has no sink. \end{theorem}
\begin{IEEEproof}
Proof of the ``if'' part of the theorem, is given by Lemma~\ref{lem:undomination}.

 For ease of notation in this proof, we use $[\Gamma]_i=\left[ \Gamma_{i1} \; \cdots \; \Gamma_{iq} \right]$ and $[A]_i=\left[ A_{i1} \; \cdots \; A_{iq} \right]$.
\par In order to prove the ``only if'' part of the theorem, we need to show that if the plant graph has a sink (i.e., if there exists $j$ such that $(s_{\mathcal{P}})_{ij}=0$ for every $i\neq j$), then there exists a control design method $\Gamma$ which dominates the deadbeat control design method. We exhibit such a strategy.

Without loss of generality, we can assume that $(s_{\mathcal{P}})_{iq}=0$ for all $i\neq q$, in which case every matrix $A$ in $\mathcal{A}(S_{\mathcal{P}})$ has the structure

\begin{equation*}
A=\matrix{cccc}{ A_{11} & \cdots & A_{1,q-1} & 0 \\ \vdots & \ddots & \vdots & \vdots \\ A_{q-1,1} & \cdots & A_{q-1,q-1} & 0 \\ A_{q1} & \cdots & A_{q,q-1} & A_{qq} }.
\end{equation*}

Define $\bar{x}_0=[\begin{array}{ccc}x_1(0) & \cdots & x_{q-1}(0) \end{array}]^T$, and let control design strategy $\Gamma$ be defined by
\begin{equation*}
\begin{split}
\matrix{cccc}{ -B_{11}^{-1}A_{11} & \cdots & -B_{11}^{-1}A_{1,q-1} & 0 \\ \vdots & \ddots & \vdots & \vdots \\ -B_{q-1,q-1}^{-1}A_{q-1,1} & \cdots & -B_{q-1,q-1}^{-1}A_{q-1,q-1} & 0 \\ K_{q1}(A,B) & \cdots & K_{q,q-1}(A,B) & K_{qq}(A,B) }
\end{split}
\end{equation*}
for all $P=(A,B,x_0) \in \mathcal{P}$, with
\begin{equation*}
\begin{split}
\bar{K}(A,B):&=\matrix{cccc}{K_{q1}(A,B) & \cdots & K_{q,q-1}(A,B) & K_{qq}(A,B)} \\
&=-(I+B_{qq}^TX_{qq}B_{qq})^{-1}B_{qq}^TX_{qq} [A]_q,
\end{split}
\end{equation*}
where $X_{qq}$ is the unique positive definite solution to the discrete algebraic Riccati equation
\begin{equation} \label{eqn:Xqq}
\begin{split}
A_{qq}^TX_{qq}B_{qq}(I+&B_{qq}^TX_{qq}B_{qq})^{-1}B_{qq}^TX_{qq}A_{qq}-A_{qq}^TX_{qq}A_{qq}+X_{qq}-I=0.
\end{split}
\end{equation}
In words, control design strategy $\Gamma$ applies the deadbeat strategy to subsystems $1$ to $q-1$ while, on subsystem $q$, it uses the same subcontroller as in the optimal controller for the plant
\begin{equation}
\label{app_sys}
\hat{x}(k+1)=\hat{A}\hat{x}(k)+\hat{B}\hat{u}(k),
\end{equation}
with cost function
\begin{equation*}
J^{(2)}_{(A,B,x_0)}(\bar{K})=\sum_{k=1}^\infty \hat{x}(k)^TQ\hat{x}(k) +\sum_{k=0}^\infty \hat{u}(k)^T\hat{u}(k),
\end{equation*}
where $Q=\mbox{diag}(0,\dots,0,I_{n_q\times n_q})$, the matrix $\hat{A}$ is defined as $[\hat{A}]_q=[A]_q$ and $[\hat{A}]_z=0$ for all $z\neq q$, and furthermore, the matrix $\hat{B}$ is defined as $\hat{B}=\mbox{diag}(0,\dots,0,B_{qq})$. Note that $\Gamma$ is indeed communication-less since $\bar{K}(A,B)$ defined above can be computed with the sole knowledge of the $q^{th}$ lower block of $A$ and $B$. Because of the structure of matrices in $\mathcal{A}(S_{\mathcal{P}})$ and this characterization of $\Gamma$, we have
\begin{equation*}
J_{(A,B,x_0)}(\Gamma(A,B))=J^{(1)}_{(A,B,x_0)}+J^{(2)}_{(A,B,x_0)}(\bar{K}(A,B)),
\end{equation*}
where $ J^{(1)}_{(A,B,x_0)}=\bar{x}_0^T \bar{A}^T \bar{B}^{-T} \bar{B}^{-1} \bar{A} \bar{x}_0, $ with
$$
\bar{A}=\matrix{ccc}{ A_{11} & \cdots & A_{1,q-1} \\ \vdots & \ddots & \vdots \\ A_{q-1,1} & \cdots & A_{q-1,q-1}},
$$
and $\bar{B}=\mbox{diag}(B_{11},\dots,B_{q-1,q-1})$ and $J^{(2)}_{(A,B,x_0)}(\bar{K}(A,B))$ is the closed-loop cost for system~(\ref{app_sys}). Since $\bar{K}(A,B)$ is the optimal controller for this cost, $J^{(2)}_{(A,B,x_0)}(\bar{K}(A,B)) = x_0^T \hat{A}^T W \hat{A} x_0$, where
\begin{equation*}
W=\mbox{diag}(0,\dots,0,X_{qq}-X_{qq}B_{qq}(I+B_{qq}^TX_{qq}B_{qq})^{-1}B_{qq}^TX_{qq}).
\end{equation*}
Using part~2 of Subsection~3.5.2 in~\cite{Handbook1996}, we have the matrix inversion identity
$$
X-XY(I+ZXY)^{-1}ZX=(X^{-1}+YZ)^{-1},
$$
which results in
\begin{equation*}
\begin{split}
W_{qq}&=X_{qq}-X_{qq}B_{qq}(I+B_{qq}^TX_{qq}B_{qq})^{-1}B_{qq}^TX_{qq} \\ &= (X_{qq}^{-1}+B_{qq}B_{qq}^T)^{-1} \\ & < B_{qq}^{-T}B_{qq}^{-1}.
\end{split}
\end{equation*}
Note that $X_{qq}^{-1}$ exists because $X_{qq}\geq I$ which follows from the discrete algebraic Riccati equation in~(\ref{eqn:Xqq}). This inequality implies that
$$
\hat{A}^TW\hat{A} < \hat{A}^T(\hat{B}^{\dag})^T\hat{B}^{\dag}\hat{A}
$$
where $\hat{B}^{\dag}=\mbox{diag}(0,\dots,0,B_{qq}^{-1})$. Thus
\begin{equation*}
\begin{split}
J_{(A,B,x_0)}(\Gamma(A,B))&=J^{(1)}_{(A,B,x_0)}+J^{(2)}_{(A,B,x_0)}(\bar{K}(A,B)) \\&< J_{(A,B,x_0)}(\Gamma^\Delta(A,B)),
\end{split}
\end{equation*}
for all $P=(A,B,x_0) \in \mathcal{P}$ such that the $q^{th}$ lower block of $A$ is not zero, unless the $J_{(A,B,x_0)}(\Gamma(A,B))=J_{(A,B,x_0)}(\Gamma^\Delta(A,B))$. Thus, control design method $\Gamma$ dominates the deadbeat control design method $\Gamma^\Delta$.
\end{IEEEproof}

\begin{remark} Consider the limited model information design problem given by the plant graph $G_\p$ in Figure~\ref{figure0}($a$), the control graph $G'_\K$ in Figure~\ref{figure0}($b'$), and the design graph $G'_\comm$ in Figure~\ref{figure0}($c'$). Theorems~\ref{tho:3} and~\ref{tho:2} show that the deadbeat control design strategy $\Gamma^\Delta$ is the best control design strategy that one can propose based on the local model of subsystems and the plant graph, because the deadbeat control design strategy is the minimizer of the competitive ratio and it is undominated.
\end{remark}

\begin{remark} It should be noted that, the proof of the ``only if'' part of the Theorem~\ref{tho:2} is constructive. We use this construction to build a control design strategy for the plant graphs with sinks in next subsection.
\end{remark}

\subsection{Second case: plant graph $G_{\mathcal{P}}$ with at least one sink}
In this section, we consider the case where plant graph $G_{\mathcal{P}}$ has $c\geq 1$ sinks. Accordingly, its adjacency matrix $S_{\mathcal{P}}$ is of the form
\begin{equation} \label{eqn:10}
S_{\mathcal{P}} = \matrix{c|c}{ (S_{\p})_{11} & 0_{(q-c)\times (c)} \\ \hline (S_{\p})_{21} & (S_{\p})_{22} },
\end{equation}
where
\begin{equation*}
(S_\mathcal{P})_{11}=\matrix{ccc}{ (s_{\p})_{11} & \cdots & (s_{\p})_{1,q-c}\\ \vdots & \ddots & \vdots \\ (s_{\p})_{q-c,1} & \cdots & (s_{\p})_{q-c,q-c} },
\end{equation*}
\begin{equation*}
(S_\mathcal{P})_{21}=\matrix{ccc}{ (s_{\p})_{q-c+1,1} & \cdots & (s_{\p})_{q-c+1,q-c}\\ \vdots & \ddots & \vdots \\ (s_{\p})_{q,1} & \cdots & (s_{\p})_{q,q-c} },
\end{equation*}
and
\begin{equation*}
(S_\mathcal{P})_{22}=\matrix{ccc}{ (s_{\p})_{q-c+1,q-c+1} & \cdots & 0\\ \vdots & \ddots & \vdots \\ 0 & \cdots & (s_{\p})_{qq} },
\end{equation*}
where we assume, without loss of generality, that the vertices are numbered such that the sinks are labeled $q-c+1,\dots,q$. With this notation, let us now introduce the control design method $\Gamma^\Theta$ defined by
\begin{equation}
\begin{split} \label{eqn:11}
\Gamma^{\Theta}(A,B)=-&\mbox{diag}(B_{11}^{-1},\dots,B_{q-c,q-c}^{-1}, W_{q-c+1}(A,B),\dots,W_q(A,B)) A
\end{split}
\end{equation}
for all $(A,B) \in \mathcal{A}(S_{\mathcal{P}}) \times \mathcal{B}(\epsilon)$, where
\begin{equation} \label{eqn:12}
W_i(A,B)=(I+B_{ii}^TX_{ii}B_{ii})^{-1}B_{ii}^TX_{ii}
\end{equation}
for all $q-c+1\leq i \leq q$ and $X_{ii}$ is the unique positive definite solution of the discrete algebraic Riccati equation
\begin{equation} \label{eqn:12.1}
\begin{split}
A_{ii}^TX_{ii}B_{ii}(I+B_{ii}^TX_{ii}&B_{ii})^{-1}B_{ii}^TX_{ii}A_{ii}-A_{ii}^TX_{ii}A_{ii}+X_{ii}-I=0.
\end{split}
\end{equation}
The control design method $\Gamma^\Theta$ applies the deadbeat strategy to every subsystem that is not a sink and, for every sink, applies the same optimal control law as if the node were decoupled from the rest of the graph. We will show that when the plant graph contains sinks, $\Gamma^{\Theta}$ has, in worst case, the same competitive ratio as the deadbeat strategy. Unlike the deadbeat strategy, it has the additional property of being undominated by communication-less methods for plants in $\mathcal{P}$ when the plant graph $G_{\mathcal{P}}$ has sinks.

\begin{lemma} \label{lem:3} Let the plant graph $G_{\mathcal{P}}$ contain no isolated node, the design graph $G_{\mathcal{C}}$ be a totally disconnected graph with self-loops, and $G_\K\supseteq G_\p$. Let $\Gamma$ be a control design strategy in $\comm$. Suppose that there exist $i$ and $j\neq i$ such that $(s_{\mathcal{P}})_{ij}\neq 0$ and that node $i$ is not a sink. The competitive ratio of $\Gamma$ is bounded only if
$$
A_{ij}+B_{ii}\Gamma_{ij}(A,B)=0, \hspace{0.1in} \textrm{for all} \; P=(A,B,x_0)\in\p.
$$ \end{lemma}
\begin{IEEEproof} See Appendix~\ref{prooflem:3}. \end{IEEEproof}

\begin{remark} Lemma~\ref{lem:3} shows that a necessary condition for a bounded competitive ratio is to decouple the nodes that are not sinks from the rest of the network. \end{remark}

 Now, we are ready to compute the competitive ratio of the newly defined control design strategy $\Gamma^\Theta$. This is done at first for the case that the control graph $G_{\mathcal{K}}$ is a complete graph.

\begin{theorem} \label{tho:5} Let the plant graph $G_{\mathcal{P}}$ contain no isolated node and at least one sink, and the control graph $G_{\mathcal{K}}$ be a complete graph. Then the competitive ratio of the communication-less design method $\Gamma^\Theta$ introduced in~(\ref{eqn:11}) is
$$
r_{\p}(\Gamma^\Theta)=\left\{ \begin{array}{ll} 1, & \mbox{if} \hspace{0.05in} (S_{\p})_{11}=0 \hspace{0.05in} \mbox{and} \hspace{0.05in} (S_{\p})_{22}=0, \\ 1+1/\epsilon^2, & \mbox{otherwise.} \end{array} \right.
$$
\end{theorem}
\begin{IEEEproof} Based on Theorem~\ref{tho:1} we know that, for every plant $P=(A,B,x_0) \in \mathcal{P}$
\begin{equation} \label{eqn:15}
\hspace{-.05in} J_{(A,B,x_0)}(K^*(A,B)) \geq \frac{\epsilon^2}{1+\epsilon^2} x_0^TA^TB^{-T}B^{-1}Ax_0,
\end{equation}
In addition, proceeding as in the proof of the ``only if'' part of the Theorem~\ref{tho:2}, we know that
\begin{equation} \label{eqn:16}
J_{(A,B,x_0)}(\Gamma^\Delta(A,B)) \geq J_{(A,B,x_0)}(\Gamma^\Theta(A,B)).
\end{equation}
Plugging equation~(\ref{eqn:16}) into equation~(\ref{eqn:15}) results in
\begin{equation*}
\frac{J_{(A,B,x_0)}(\Gamma^\Theta(A,B))}{J_{(A,B,x_0)}(K^*(A,B))} \leq 1+ \frac{1}{\epsilon^2}, \; \forall P=(A,B,x_0) \in \mathcal{P}.
\end{equation*}
As a result, $r_{\mathcal{P}} (\Gamma^{\Theta}) \leq 1+ 1/\epsilon^2$. To show that this upper-bound is tight, we now exhibit plants for which it is attained. We use a different construction depending on matrices $(S_{\p})_{11}$ and $(S_{\p})_{22}$. If $(S_{\p})_{11}\neq 0$, two situations can occur.
\\ \textit{Case 1: $(S_{\p})_{11}\neq 0$ and it is not diagonal.} There exist $1\leq i \neq j\leq q-c$ such that $(s_{\mathcal{P}})_{ij}\neq 0$. In this case, choose indices $i_1\in \mathcal{I}_i$ and $j_1\in \mathcal{I}_j$ and define $A=e_{i_1}e_{j_1}^T$ and $B=\epsilon I$. Then, for $x_0=e_{j_1}$, we find that
\begin{equation*}
\frac{J_{(A,B,x_0)}(\Gamma^\Theta(A,B))}{J_{(A,B,x_0)}(K^*(A,B))} = \frac{1/\epsilon^2}{1/(1+\epsilon^2)} = 1+ \frac{1}{\epsilon^2}
\end{equation*}
because the control design $\Gamma^\Theta$ acts like the deadbeat control design method on this plant.
\\ \textit{Case 2: $(S_{\p})_{11}\neq 0$ and it is diagonal.} There exists $1\leq i\leq q-c$ such that $(s_{\mathcal{P}})_{ii}\neq 0$. Pick an index $i_1\in\mathcal{I}_i$. In that case, consider $A(r)=re_{i_1}e_{i_1}^T$ and $B=\epsilon I$. For $x_0=e_{i_1}$, the optimal cost is
\begin{equation*}
\begin{split}
&J_{(A(r),B,x_0)}(K^*(A(r),B))= \frac{\sqrt{r^4 + 2r^2\epsilon^2 - 2r^2 + \epsilon^4 + 2\epsilon^2 + 1}+r^2 - \epsilon^2 - 1}{2\epsilon^2},
\end{split}
\end{equation*}
which results in
\begin{equation*}
\lim_{r\rightarrow 0} \frac{J_{(A,B,x_0)}(\Gamma^\Theta(A,B))}{J_{(A,B,x_0)}(K^*(A,B))} = 1+ \frac{1}{\epsilon^2}.
\end{equation*}
\\ Now suppose that $(S_{\p})_{11}=0$. Again, two different situations can occur.
\\ \textit{Case 3: $(S_{\p})_{11}=0$ and $(S_{\p})_{22}\neq 0$.} There exists $q-c+1\leq i\leq q$ such that $(s_\p)_{ii}\neq 0$. From the assumption that the plant graph contains no isolated node, we know that there must exist $1 \leq j \leq q-c$ such that $(s_{\mathcal{P}})_{ij} \neq 0$. Accordingly, let us pick $i_1 \in \mathcal{I}_i$ and $j_1 \in \mathcal{I}_j$ and consider the $2$-parameter family of matrices $A(r,s)$ in $\mathcal{A}(S_{\mathcal{P}})$ with all entries equal to zero except $a_{i_1 i_1}$, which is equal to $r$, and $a_{i_1 j_1}$, which is equal to $s$. Let $B=\epsilon I$. For any initial condition $x_0$, the corresponding closed-loop performance is
\begin{equation*}
J_{(A(r,s),B,x_0)}(\Gamma^\Theta(A(r,s),B))=\beta_\Theta x_0^Ta(r,s)a(r,s)^Tx_0,
\end{equation*}
where
we have let $a(r,s)=A(r,s)_{i_1}^T$ and $\beta_\Theta$ is
\begin{equation*}
\beta_\Theta=\frac{\sqrt{r^4 + 2r^2\epsilon^2 - 2ar^2 + \epsilon^4 + 2\epsilon^2 + 1}+r^2 - \epsilon^2 - 1}{2\epsilon^2 r^2}.
\end{equation*}
Besides, the optimal closed-loop performance can be computed as
\begin{equation*}
J_{(A(r,s),B,x_0)}(K^*(A(r,s),B))=\beta_{K^*} x_0^Ta(r,s)a(r,s)^Tx_0,
\end{equation*}
where $\beta_{K^*}$ is
\begin{equation*}
\beta_{K^*}=\frac{\epsilon^2s^2 +r^2(1+\epsilon^2)-(\epsilon^2+1)^2 +\sqrt{c_+c_-}}{2\epsilon^2(\epsilon^2 + 1)(s^2 + r^2)},
\end{equation*}
\begin{equation*}
c_\pm=(\epsilon^2s^2+ (r^2\pm 2r)(\epsilon^2+1) + (\epsilon^2+1)^2).
\end{equation*}
Then,
\begin{equation*}
\begin{split}
r_{\p}(\Gamma^{\Theta}) &\geq \lim_{r\rightarrow \infty, \frac{s}{r} \rightarrow \infty} \frac{J_{(A(r,s),B,x_0)}(\Gamma^\Theta(A(r,s),B))}{J_{(A(r,s),B,x_0)}(K^*(A(r,s),B))} \\&=1+\frac{1}{\epsilon^2}
\end{split}
\end{equation*}
\\ \textit{Case 4: $(S_{\p})_{11}=0$ and $(S_{\p})_{22}=0$.} Then, every matrix $A \in \mathcal{A}(S_{\mathcal{P}})$ has the form
$
\left[
\begin{array}{c|c}
0 & 0\\
\hline
* & 0
\end{array}
\right]
$
and, in particular, is nilpotent of degree 2; i.e., $A^2=0$. In this case, the Riccati equation yielding the optimal control gain $K^*(A,B)$ can be readily solved, and we find that $K^*(A,B)=-(I+B^TB)^{-1}B^TA$ for all $(A,B)$. As a result, $K^*(A,B)=\Gamma^{\Theta}(A,B)$ for all plant $P=(A,B,x_0) \in \mathcal{P}$ (since $W_i(A,B) = (I+B_{ii}^TB_{ii})^{-1}B_{ii}^T$ for all $q-c+1 \leq i \leq q$), which implies that the competitive ratio of $\Gamma^{\Theta}$ against plants in $\mathcal{P}$ is equal to one.
\end{IEEEproof}

In Theorem~\ref{tho:5}, the control graph $G_\K$ is assumed to be a complete graph. We needed this assumption to calculate the cost of the optimal control design strategy $K^*(P)$ when $(S_\p)_{11}=0$ and $(S_\p)_{22}\neq 0$ which is not an easy task when the control graph $G_\K$ is incomplete. However, more can be said if $(S_\p)_{11}\neq 0$.

\begin{corollary} Let the plant graph $G_{\mathcal{P}}$ contain no isolated node and at least one sink and $G_\K\supseteq G_\p$. Then
$$
r_{\p}(\Gamma^\Theta)=\left\{ \begin{array}{ll} 1, & \mbox{if } (S_{\p})_{11}=0 \mbox{ and } (S_{\p})_{22}=0, \\ 1+1/\epsilon^2, & \mbox{if } (S_{\p})_{11}\neq 0. \end{array} \right.
$$
\end{corollary}
\begin{IEEEproof} According to Theorem~\ref{tho:5}, for $(S_\p)_{11}\neq 0$, we get
\begin{equation*}
\begin{split}
r_\p(\Gamma^\Theta) &= \sup_{P\in\p}\frac{J_{(A,B,x_0)}(\Gamma^\Theta(A,B))}{J_{(A,B,x_0)}(K^*(P))} \\&\leq \sup_{P\in\p}\frac{J_{(A,B,x_0)}(\Gamma^\Theta(A,B))}{J_{(A,B,x_0)}(K^*_C(A,B))} = 1+ \frac{1}{\epsilon^2}.
\end{split}
\end{equation*}
\\ \textit{Case 1: $(S_{\p})_{11}\neq 0$ and it is not diagonal.} For the special plant introduced in Case~1 in the proof of Theorem~\ref{tho:5}, we have
$J_{(A,B,e_{j_1})}(K^*_C(A,B))=\linebreak[4]J_{(A,B,e_{j_1})}(K^*(A,B,e_{j_1}))$ since $A=e_{i_1}e_{j_1}^T$ is a nilpotent matrix. The rest of the proof is similar to Case~1 in the proof of Theorem~\ref{tho:5}.
\\ \textit{Case 2: $(S_{\p})_{11}\neq 0$ and it is diagonal.} Note that, for the special plant introduced Case~2 in the proof of Theorem~\ref{tho:5}, we have
\begin{equation*}
\begin{split}
&K^*_C(A,B)=-\frac{\sqrt{r^4 + 2r^2\epsilon^2 - 2r^2 + \epsilon^4 + 2\epsilon^2 + 1}+r^2 - \epsilon^2 - 1}{2\epsilon r^2} A
\end{split}
\end{equation*}
which shows $K^*_C(A,B)\in\K(S_\K)$ and similar to the proof of Theorem~\ref{tho:1}, we get $J_{(A,B,e_{i_1})}(K^*_C(A,B))=J_{(A,B,e_{i_1})}(K^*(A,B,e_{i_1}))$. The rest of the proof is similar to Case~2 in the proof of Theorem~\ref{tho:5}.
\\ \textit{Case 3: $(S_{\p})_{11}=0$ and $(S_{\p})_{22}=0$.} Then, every $A \in \mathcal{A}(S_{\mathcal{P}})$ is nilpotent matrix which results in $J_P(K^*(P))=J_P(K^*_C(A,B))$. The rest of the proof is similar to Case~4 in the proof of Theorem~\ref{tho:5}.
\end{IEEEproof}

Now that we have computed the competitive ratio of the control design strategy $\Gamma^{\Theta}$ in the presence of sinks, we present a theorem to show that the competitive ratio of any other communication-less control design strategy is lower-bounded by the competitive ratio of $\Gamma^{\Theta}$ when the control graph $G_\mathcal{K}$ is a complete graph. Therefore, the control design strategy $\Gamma^{\Theta}$ is a minimizer of the competitive ratio over the set of limited model information control design strategies.

\begin{theorem} \label{tho:6} Let the plant graph $G_{\mathcal{P}}$ contain no isolated node and at least one sink, the control graph $G_{\mathcal{K}}$ be a complete graph, and the design graph $G_{\mathcal{C}}$ be a totally disconnected graph with self-loops. Then the competitive ratio of any control design strategy $\Gamma\in\comm$ satisfies
$$
r_{\p}(\Gamma)\geq 1+1/\epsilon^2,
$$
if either $(S_{\p})_{11}$ is not diagonal or $(S_{\p})_{22}\neq 0$. \end{theorem}
\begin{IEEEproof} \textit{Case 1: $(S_{\p})_{11}\neq 0$ and it is not diagonal.} Then, there exist $1\leq i,j \leq q-c$ and $i\neq j$ such that $(s_{\mathcal{P}})_{ij}\neq 0$. Choose indices $i_1 \in \mathcal{I}_i$ and $j_1 \in \mathcal{I}_j$ and consider the matrix $A$ defined by $A=e_{i_1}e_{j_1}^T$ and $B=\epsilon I$. From Lemma~\ref{lem:3}, we know that a communication-less method $\Gamma$ has a bounded competitive ratio only if $\Gamma(A,B)=-B^{-1}A$ (because node $i$ is a part of $(S_{\p})_{11}$ and it is not a sink). Therefore
\begin{equation*}
r_{\p}(\Gamma) \geq \frac{J_{(A,B,e_{j_1})}(\Gamma(A,B))}{J_{(A,B,e_{j_1})}(K^*(A,B))}=1+\frac{1}{\epsilon^2}
\end{equation*}
for any such method. \\
\textit{Case 2: $(S_{\p})_{22}\neq 0$.} There thus exists $q-c+1\leq i\leq q$ such that $(s_{\mathcal{P}})_{ii}\neq 0$. Note that, there exists $1\leq j\leq q-c$ such that $(s_{\mathcal{P}})_{ij}\neq 0$, since there is no isolated node in the plant graph. Choose indices $i_1 \in \mathcal{I}_i$ and $j_1 \in \mathcal{I}_j$. Consider $A$ defined as $A=re_{i_1}e_{j_1}^T+se_{i_1}e_{i_1}^T$ and $B=\epsilon I$. As indicated in the proof of Theorem~\ref{tho:5}, control design strategy $\Gamma^\Theta$ yields the globally optimal controller with limited model information for plants in this family. Hence, we know that $r_\p(\Gamma) \geq 1 +1/\epsilon^2$ for every communication-less strategy $\Gamma$.
\end{IEEEproof}

 In Theorem~\ref{tho:6}, we assume the control graph $G_\K$ is a complete graph. In the next corollary, we generalize this result to the case where $G_\K$ is a supergraph of $G_\p$ when $(S_{\p})_{11}$ is not diagonal.

\begin{corollary} Let the plant graph $G_{\mathcal{P}}$ contain no isolated node and at least one sink, the design graph $G_{\mathcal{C}}$ be a totally disconnected graph with self-loops, and $G_\K \supseteq G_\p$. Then the competitive ratio of any control design strategy $\Gamma\in\comm$ satisfies
$$
r_{\p}(\Gamma)\geq 1+1/\epsilon^2,
$$
if $(S_{\p})_{11}$ is not diagonal.
\end{corollary}
\begin{IEEEproof} Considering that for the nilpotent matrix $A=e_{i_1}e_{j_1}^T$, we get \linebreak[4] $J_{(A,B,e_{j_1})}(K^*(A,B,e_{j_1}))=J_{(A,B,e_{j_1})}(K^*_C(A,B))$, the rest of the proof is similar to Case~1 in the proof of Theorem~\ref{tho:6}. \end{IEEEproof}

\begin{remark} Combining Theorems~\ref{tho:5} and~\ref{tho:6} implies that if either $(S_{\p})_{11}$ is not diagonal or $(S_{\p})_{22}\neq 0$, control design method $\Gamma^{\Theta}$ exhibits the same competitive ratio as the deadbeat control strategy, which is the smallest ratio achievable by a communication-less control method. Therefore, it is a solution to problem~(\ref{eqn:0}). Furthermore, if $(S_{\p})_{11}$ and $(S_{\p})_{22}$ are both zero, then $\Gamma^\Theta$ is equal to $K^*$, which shows that $\Gamma^\Theta$ is a solution to problem~(\ref{eqn:0}), in this case too. \end{remark}

\begin{remark} The case where $(S_{\p})_{11}$ is diagonal and $(S_{\p})_{22}=0$ is still open. \end{remark}

The next theorem shows that $\Gamma^{\Theta}$ is a more desirable control design method than the deadbeat control design strategy when the plant graph $G_{\mathcal{P}}$ has sinks, since it is then undominated by communication-less design methods.

\begin{theorem} \label{tho:4} Let the plant graph $G_{\mathcal{P}}$ contain no isolated node and at least one sink, the design graph $G_{\mathcal{C}}$ be a totally disconnected graph with self-loops, and $G_{\mathcal{K}}\supseteq G_\p$. The control design method $\Gamma^\Theta$ is undominated by any control design method $\Gamma\in\comm$. \end{theorem}
\begin{IEEEproof} See Appendix~\ref{prooftho:4}. \end{IEEEproof}

\begin{remark} Consider the limited model information design problem given by the plant graph $G'_\p$ in Figure~\ref{figure0}($a'$), the control graph $G'_\K$ in Figure~\ref{figure0}($b'$), and the design graph $G'_\comm$ in Figure~\ref{figure0}($c'$). Theorems~\ref{tho:5},~\ref{tho:6}, and~\ref{tho:4} together show that, the control design strategy $\Gamma^\Theta$ is the best control design strategy that one can propose based on the local model information and the plant graph, because the control design strategy $\Gamma^\Theta$ is a minimizer of the competitive ratio and it is undominated.
\end{remark}

\begin{remark} For general weight matrices $Q$ and $R$ appearing in the performance cost, the competitive ratio of both the deadbeat control design strategy $\Gamma^\Delta$ and the control design strategy $\Gamma^\Theta$ is $1+\bar{\sigma}(R)/(\underline{\sigma}(Q)\epsilon^2)$. In particular, the competitive ratio has a limit equal to one as $\bar{\sigma}(R)/\underbar{\sigma}(Q)$ goes to zero. We thus recover the well-known observation (e.g.,~\cite{O'Reilly1981363}) that, for discrete-time linear time-invariant systems, the optimal linear quadratic regulator approaches the deadbeat controller in the limit of ``cheap control''.
\end{remark}

\section{Design Graph Influence on Achievable Performance} \label{sec:3}
In the previous section, we have shown that communicat-ion-less control design methods (i.e., $G_{\mathcal{C}}$ is totally disconnected with self-loops) have intrinsic performance limitations, and we have characterized minimal elements for both the competitive ratio and domination metrics. A natural question is ``given plant graph $G_{\mathcal{P}}$, which design graph $G_{\mathcal{C}}$ is necessary to ensure the existence of $\Gamma \in \mathcal{C}$ with better competitive ratio than $\Gamma^{\Delta}$ and $\Gamma^{\Theta}$ ?''. We tackle this question in this section.

\begin{theorem} \label{lem:4} Let the plant graph $G_\p$ and the design graph $G_\comm$ be given and $G_\mathcal{K}\supseteq G_\p$. If one of the following conditions is satisfied then $r_{\p}(\Gamma)\geq 1+1/\epsilon^2$ for all $\Gamma \in \mathcal{C}$:
\begin{enumerate}
\item[ (a) ] $G_{\mathcal{P}}$ contains the path $k\rightarrow i \rightarrow j$ with distinct nodes $i$, $j$, and $k$ while $(j,i)\notin E_\comm$.
\item[ (b) ] There exist $i\neq j$ such that $n_i \geq 2$ and $(i,j) \in E_{\mathcal{P}}$ while $(j,i)\notin E_\comm$.
\end{enumerate}
\end{theorem}

\begin{IEEEproof}
We prove the case when condition~(a) holds. The proof for condition~(b) is similar.

Let $i$, $j$, and $k$ be three distinct nodes such that $(s_{\mathcal{P}})_{ik}\neq 0$ and $(s_{\mathcal{P}})_{ji}\neq 0$ (i.e., the path $k\rightarrow i \rightarrow j$ is contained in the plant graph $G_{\mathcal{P}}$). Let us pick $i_1\in\mathcal{I}_i$, $j_1\in \mathcal{I}_j$ and $k_1\in \mathcal{I}_k$ and consider the 2-parameter family of matrices $A(r,s)$ in $\mathcal{A}(S_\p)$ with all entries equal to zero except $a_{i_1k_1}$, which is equal to $r$, and $a_{j_1i_1}$, which is equal to $s$. Let $B=\epsilon I$ and let $\Gamma\in \comm$ be a limited model information with design graph $G_{\mathcal{C}}$. For $x_0=e_{k_1}$, we have
\begin{equation*}
\begin{split}
J_{(A(r,s),B,e_{k_1})}(\Gamma(A(r,s)&,B)) \geq (r+\epsilon \gamma_{i_1k_1}(A,B))^2 [\gamma_{j_1i_1}^2+(s+\epsilon \gamma_{j_1i_1}(A,B))^2]
\end{split}
\end{equation*}
where $\gamma_{i_1k_1}$ cannot be a function of $s$ because $(j,i)\notin E_\comm$. Note that, irrespective of the choice of $\gamma_{j_1i_1}(A,B)$, we have
\begin{equation*}
\begin{split}
J_{(A(r,s),B,e_{k_1})}(\Gamma(A(r,s)&,B)) \geq \frac{(r+\epsilon \gamma_{i_1k_1}(A,B))^2s^2}{1+\epsilon^2}.
\end{split}
\end{equation*}
The cost of the deadbeat control design on this plant satisfies
\begin{equation*}
\begin{split}
J_{(A(r,s),B,e_{k_1})}(\Gamma^\Delta(A(r,s)&,B)) = r^2/\epsilon^2,
\end{split}
\end{equation*}
and thus
\begin{equation} \label{inequalityusingdeadbeat}
\begin{split}
r_{\p}(\Gamma) &=\sup_{P\in \p} \frac{J_{P}(\Gamma(A,B))}{J_{P}(K^*(P))} \\ &=\sup_{P\in \p} \left[ \frac{J_{P}(\Gamma(A,B))}{J_{P}(\Gamma^\Delta(A,B))} \frac{J_{P}(\Gamma^\Delta(A,B))}{J_{P}(K^*(P))} \right] \\ & \geq \sup_{P\in \p} \frac{J_{P}(\Gamma(A,B))}{J_{P}(\Gamma^\Delta(A,B))},\\
& \geq \lim_{s\rightarrow \infty} \frac{\epsilon^2(r+\epsilon \gamma_{i_1k_1}(A,B))^2s^2}{(1+\epsilon^2) r^2}.
\end{split}
\end{equation}
This shows that $r_\p(\Gamma)$ is unbounded unless $r+\epsilon \gamma_{i_1k_1}(A(r,s),B)=0$ for all $r,s$. Now consider the 1-parameter family of matrices $\bar{A}(r)$ with all entries equal to zero except $a_{i_1k_1}$, which is equal to $r$. Because of $(j,i)\notin E_\comm$, we know that $\Gamma_z(\bar{A}(r),B)=\Gamma_z(A(r,s),B)$ for all $z\in \mathcal{I}_i$. Thus
\begin{equation*}
J_{(\bar{A}(r),B,e_{k_1})}(\Gamma(\bar{A}(r),B)) \geq r^2/\epsilon^2.
\end{equation*}
On the other hand, similar to the proof of Theorem~\ref{tho:1}, we can compute the optimal controller for systems in this $1-$parameter family and find
\begin{equation*}
\begin{split}
J_{(\bar{A}(r),B,e_{k_1})}(K^*(\bar{A}(r),B,e_{k_1}))&=J_{(\bar{A}(r),B,e_{k_1})}(K^*_C(\bar{A}(r),B))\\&=r^2/(1+\epsilon^2),
\end{split}
\end{equation*}
As a result, we get
\begin{equation*}
r_\p(\Gamma) \geq \frac{r^2/\epsilon^2}{r^2/(1+\epsilon^2)}=1+\frac{1}{\epsilon^2},
\end{equation*}
which concludes the proof for this case.
\end{IEEEproof}

\begin{remark} Consider the limited model information design problem given by the plant graph $G_\p$ in Figure~\ref{figure0}($a$), the control graph $G'_\K$ in Figure~\ref{figure0}($b'$), and the design graph $G_\comm$ in Figure~\ref{figure0}($c$). Theorem~\ref{lem:4} shows that, because the plant graph $G_\mathcal{P}$ contains the path $3\rightarrow2\rightarrow1$ but the design graph $G_\mathcal{C}$ does not contain $1\rightarrow 2$, the competitive ratio of any control design strategy $\Gamma\in\comm$ would be greater than or equal to $1+1/\epsilon^2$.
\end{remark}

\begin{corollary}\label{tho:7}
Let both the plant graph $G_{\mathcal{P}}$ and the control graph $G_\mathcal{K}$ be complete graphs. If the design graph $G_{\mathcal{C}}$ is not equal to $G_{\mathcal{P}}$, then $r_{\mathcal{P}}(\Gamma) \geq 1 + 1/\epsilon^2$ for all $\Gamma \in \mathcal{C}$.
\end{corollary}
\begin{IEEEproof} The proof is a direct application of Theorem~\ref{lem:4} with condition~(a) fulfilled. \end{IEEEproof}

\begin{remark} Corollary~\ref{tho:7} shows that, when $G_{\mathcal{P}}$ is a complete graph, achieving a better competitive ratio than the deadbeat design strategy requires each subsystem to have full knowledge of the plant model when constructing each subcontroller. \end{remark}

\section{Extensions to Under-Actuated Sinks } \label{subsec:6.1}
In the previous sections, we gave an explicit solution to the problem in~(\ref{eqn:0}) under the assumption that all the subsystems are fully-actuated; i.e., all the matrices $B\in\B(\epsilon)$ are square invertible matrices. Note that this assumption stems from the fact that the subsystems that are not sinks in the plant graph are required to decouple themselves from the rest of the plant to avoid influencing highly sensitive (and potentially hard to control) subsystems in order to keep the competitive ratio finite (see Lemma~\ref{lem:3}). Therefore, we assume these subsystems are fully-actuated to easily decouple them from the rest of the system. As a future direction for improvement, one can try to replace this assumption with other conditions (e.g.,~geometric conditions) to ensure that the subsystems can decouple themselves. From the same argument, it should be expected that the assumption of a square invertible B-matrix is dispensable for sink nodes. In this section, we briefly discuss an extension of our results to the slightly more general, but still restricted, class of plants whose sinks are under-actuated.

Consider the limited model information control design problem given with the plant graph $G_\p$, the control graph $G_\K$, and the design graph $G_\comm$ given in Figure~\ref{figure_inv_B}. The state space representation of the system is given as
$$
\matrix{c}{\underline{x}_1(k+1)\\ \underline{x}_2(k+1)}=A \matrix{c}{\underline{x}_1(k)\\ \underline{x}_2(k)}+B \matrix{c}{\underline{u}_1(k) \\ \underline{u}_2(k)},
$$
where
$$
A=\matrix{cc}{A_{11} & 0 \\ A_{21} & A_{22}}, \; B=\matrix{cc}{B_{11} & 0 \\ 0 & B_{22}},
$$
with $\underline{x}_1(k)\in\mathbb{R}^{n_1}$, $\underline{x}_2(k)\in\mathbb{R}^{n_2}$, $\underline{u}_1(k)\in\mathbb{R}^{n_1}$, and $\underline{u}_2(k)\in\mathbb{R}^{m_2}$ for some given integers $n_1\geq 1$, $n_2>m_2\geq 1$. Thus, for the second subsystem the matrix $B_{22}\in \mathbb{R}^{n_2\times m_2}$ is a non-square matrix, and as a result the second subsystem is an under-actuated subsystem. Let us assume that the matrices $A_{21}$, $A_{22}$, $B_{22}$ satisfy the ``matching condition''; i.e., the pair $(A_{22},B_{22})$ is controllable and $\mbox{span}(A_{21})\subseteq \mbox{span}(B_{22})$~\cite{Siljak1991}. Besides, assume that for all matrices $B$, we have $\underline{\sigma}(B) \geq \epsilon$ for some $\epsilon>0$. For this case, we have
$$
\Gamma^\Theta(A,B)=-\diag(B_{11}^{-1},W_2(A_{22},B_{22}))A,
$$
where $W_2(A_{22},B_{22})$ is defined in~(\ref{eqn:12}). Note that we do not require the matrix $B_{22}$ to be square invertible. Under some additional conditions and following a similar approach as above, it can be shown that the control design strategy $\Gamma^\Theta$ becomes an undominated minimizer of the competitive ratio over the set of limited model information control design strategies. This result can be generalized to cases with higher number of subsystems as long as the sinks in the plant graph $G_\mathcal{P}$ are the only under-actuated subsystems~\cite{Farokhi2011}.

\jpfig{0.6}{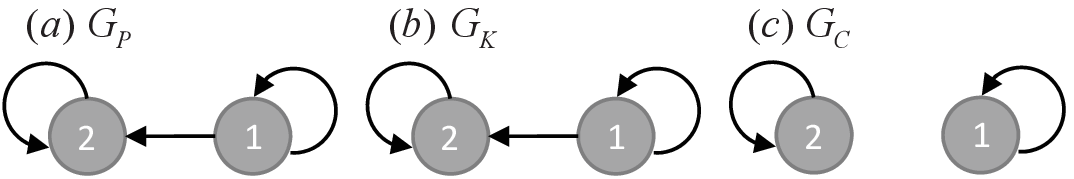}{figure_inv_B}{Plant graph $G_\mathcal{P}$, control graph $G_\mathcal{K}$, and design graph $G_\mathcal{C}$ used to illustrate an extension to under-actuated systems.}

\section{Conclusion} \label{sec:7}
We presented a framework for the study of control design under limited model information, and investigated the connection between the quality of controllers produced by a design method and the amount of plant model information available to it. We showed that the best performance achievable by a limited model information control design method crucially depends on the structure of the plant graph and, thus, that giving the designer access to this graph, even without a detailed model of all plant subsystems, results in superior design, in the sense of domination. Possible future work will focus on extending the present framework to dynamic controllers and/or where disturbances are present.

\bibliographystyle{ieeetr}
\bibliography{compile_new}

\begin{thebibliography}{10}

\bibitem{Farokhi10}
F.~Farokhi, C.~Langbort, and K.~H. Johansson, ``Control design with limited
  model information,'' in {\em American Control Conference, Proceedings of
  the}, pp.~4697 -- 4704, 2011.

\bibitem{Giulietti2000}
F.~Giulietti, L.~Pollini, and M.~Innocenti, ``Autonomous formation flight,''
  {\em Control Systems Magazine, IEEE}, vol.~20, no.~6, pp.~34 -- 44, 2000.

\bibitem{Kapilal1999}
V.~Kapila, A.~G. Sparks, J.~M. Buffington, and Q.~Yan, ``Spacecraft formation
  flying: dynamics and control,'' in {\em American Control Conference,
  Proceedings of the}, vol.~6, pp.~4137 -- 4141, 1999.

\bibitem{Swaroop1999}
D.~Swaroop and J.~K. Hedrick, ``Constant spacing strategies for platooning in
  automated highway systems,'' {\em Journal of Dynamic Systems, Measurement,
  and Control}, vol.~121, no.~3, pp.~462 -- 470, 1999.

\bibitem{Negenborn2010}
R.~R. Negenborn, Z.~Lukszo, and H.~Hellendoorn, eds., {\em Intelligent
  Infrastructures}, vol.~42.
\newblock Springer, 2010.

\bibitem{Joshi1989}
S.~M. Joshi, {\em Control of large flexible space structures}.
\newblock Springer-Verlag, 1989.

\bibitem{Dunbar2007}
W.~B. Dunbar, ``Distributed receding horizon control of dynamically coupled
  nonlinear systems,'' {\em Automatic Control, IEEE Transactions on}, vol.~52,
  no.~7, pp.~1249 -- 1263, 2007.

\bibitem{Levine1971}
W.~Levine, T.~Johnson, and M.~Athans, ``Optimal limited state variable feedback
  controllers for linear systems,'' {\em Automatic Control, IEEE Transactions
  on}, vol.~16, no.~6, pp.~785 -- 793, 1971.

\bibitem{Wenk1980}
C.~Wenk and C.~Knapp, ``Parameter optimization in linear systems with
  arbitrarily constrained controller structure,'' {\em Automatic Control, IEEE
  Transactions on}, vol.~25, no.~3, pp.~496 -- 500, 1980.

\bibitem{Soderstrom1978}
T.~S\"{o}derstr\"{o}m, ``On some algorithms for design of optimal constrained
  regulators,'' {\em Automatic Control, IEEE Transactions on}, vol.~23, no.~6,
  pp.~1100 -- 1101, 1978.

\bibitem{Shih-HoWang1973}
S.-H. Wang and E.~Davison, ``On the stabilization of decentralized control
  systems,'' {\em Automatic Control, IEEE Transactions on}, vol.~18, no.~5,
  pp.~473 -- 478, 1973.

\bibitem{Paganini1999}
G.~Ayres~de Castro and F.~Paganini, ``Control of distributed arrays with
  recursive information flow: some case studies,'' in {\em Decision and
  Control, Proceedings of the 38th IEEE Conference on}, vol.~1, pp.~191 -- 196,
  1999.

\bibitem{Bamieh2002}
B.~Bamieh, F.~Paganini, and M.~A. Dahleh, ``Distributed control of spatially
  invariant systems,'' {\em Automatic Control, IEEE Transactions on}, vol.~47,
  no.~7, pp.~1091 -- 1107, 2002.

\bibitem{Lall2003}
B.-D. Chen and S.~Lall, ``Dissipation inequalities for distributed systems on
  graphs,'' in {\em Decision and Control, Proceedings of the 42nd IEEE
  Conference on}, vol.~3, pp.~3084 -- 3090, 2003.

\bibitem{Hu1994}
Z.~Hu, ``Decentralized stabilization of large scale interconnected systems with
  delays,'' {\em Automatic Control, IEEE Transactions on}, vol.~39, no.~1,
  pp.~180 -- 182, 1994.

\bibitem{Scorletti2001}
G.~Scorletti and G.~Duc, ``An \textsc{LMI} approach to decentralized control,''
  {\em International Journal of Control}, vol.~74, no.~3, pp.~211 -- 224, 2001.

\bibitem{Rotkowitz2006}
M.~Rotkowitz and S.~Lall, ``A characterization of convex problems in
  decentralized control,'' {\em Automatic Control, IEEE Transactions on},
  vol.~51, no.~2, pp.~274 -- 286, 2006.

\bibitem{Voulgaris200351}
P.~G. Voulgaris, ``Optimal control of systems with delayed observation sharing
  patterns via input-output methods,'' {\em Systems \& Control Letters},
  vol.~50, no.~1, pp.~51 -- 64, 2003.

\bibitem{Doyle1982}
J.~C. Doyle, ``Analysis of feedback systems with structured uncertainties,''
  {\em Control Theory and Applications, IEE Proceedings D}, vol.~129, no.~6,
  pp.~242 -- 250, 1982.

\bibitem{Zames1981}
G.~Zames, ``Feedback and optimal sensitivity: Model reference transformations,
  multiplicative seminorms, and approximate inverses,'' {\em Automatic Control,
  IEEE Transactions on}, vol.~26, no.~2, pp.~301 -- 320, 1981.

\bibitem{Ball1987}
J.~A. Ball and N.~Cohen, ``Sensitivity minimization in an \textsc{H}$_\infty$
  norm: parametrization of all suboptimal solutions,'' {\em International
  Journal of Control}, vol.~46, no.~3, pp.~785 -- 816, 1987.

\bibitem{zhou1998}
K.~Zhou and J.~C. Doyle, {\em Essentials of robust control}.
\newblock Prentice Hall, 1998.

\bibitem{Ando1963}
A.~Ando and F.~M. Fisher, ``Near-decomposability, partition and aggregation,
  and the relevance of stability discussions,'' {\em International Economic
  Review}, vol.~4, no.~1, pp.~53 -- 67, 1963.

\bibitem{Sethi1998}
S.~Sethi and Q.~Zhang, ``Near optimization of dynamic systems by decomposition
  and aggregation,'' {\em Journal of Optimization Theory and Applications},
  vol.~99, no.~1, pp.~1 -- 22, 1998.

\bibitem{Gajtsgori1979}
V.~G. Gajtsgori and A.~A. Pervozvanski, ``Perturbation method in the optimal
  control problems,'' {\em Systems Science}, vol.~5, pp.~91 -- 102, 1979.

\bibitem{Sezer1986}
M.~E. Sezer and D.~D. \v{S}iljak, ``Nested $\epsilon$-decompositions and
  clustering of complex systems,'' {\em Automatica}, vol.~22, no.~3, pp.~321 --
  331, 1986.

\bibitem{langbort10}
C.~Langbort and J.-C. Delvenne, ``Distributed design methods for linear
  quadratic control and their limitations,'' {\em Automatic Control, IEEE
  Transactions on}, vol.~55, no.~9, pp.~2085 -- 2093, 2010.

\bibitem{Baughman1997}
M.~Baughman, S.~Siddiqi, and J.~Zarnikau, ``Advanced pricing in electrical
  systems. i. theory,'' {\em Power Systems, IEEE Transactions on}, vol.~12,
  no.~1, pp.~489 -- 495, 1997.

\bibitem{Berger1989}
A.~Berger and F.~Schweppe, ``Real time pricing to assist in load frequency
  control,'' {\em Power Systems, IEEE Transactions on}, vol.~4, no.~3, pp.~920
  -- 926, 1989.

\bibitem{ChaoPeck1996}
H.-P. Chao and S.~Peck, ``A market mechanism for electric power transmission,''
  {\em Journal of Regulatory Economics}, vol.~10, pp.~25--59, 1996.

\bibitem{Botterud2005}
A.~Botterud, M.~Ilic, and I.~Wangensteen, ``Optimal investments in power
  generation under centralized and decentralized decision making,'' {\em Power
  Systems, IEEE Transactions on}, vol.~20, no.~1, pp.~254 -- 263, 2005.

\bibitem{Braun2003229}
M.~Braun, D.~Rivera, M.~Flores, W.~Carlyle, and K.~Kempf, ``A model predictive
  control framework for robust management of multi-product, multi-echelon
  demand networks,'' {\em Annual Reviews in Control}, vol.~27, no.~2, pp.~229
  -- 245, 2003.

\bibitem{ChenStankovic2005}
X.-B. Chen and S.~S. Stankovic, ``Decomposition and decentralized control of
  systems with multi-overlapping structure,'' {\em Automatica}, vol.~41,
  no.~10, pp.~1765 -- 1772, 2005.

\bibitem{Zeynelgil2002}
H.~Zeynelgil, A.~Demiroren, and N.~Sengor, ``The application of ann technique
  to automatic generation control for multi-area power system,'' {\em
  International Journal of Electrical Power \& Energy Systems}, vol.~24, no.~5,
  pp.~345 -- 354, 2002.

\bibitem{PhysRevLett751226}
T.~Vicsek, A.~Czir\'{o}k, E.~Ben-Jacob, I.~Cohen, and O.~Shochet, ``Novel type
  of phase transition in a system of self-driven particles,'' {\em Physical
  Review Letters}, vol.~75, pp.~1226--1229, 1995.

\bibitem{Jadbabaie1205192}
A.~Jadbabaie, J.~Lin, and A.~S. Morse, ``Coordination of groups of mobile
  autonomous agents using nearest neighbor rules,'' {\em Automatic Control,
  IEEE Transactions on}, vol.~48, no.~6, pp.~988 -- 1001, 2003.

\bibitem{Emami-Naeini1982}
A.~Emami-Naeini and G.~Franklin, ``Deadbeat control and tracking of
  discrete-time systems,'' {\em Automatic Control, IEEE Transactions on},
  vol.~27, no.~1, pp.~176 -- 181, 1982.

\bibitem{Levine1970}
W.~Levine and M.~Athans, ``On the determination of the optimal constant output
  feedback gains for linear multivariable systems,'' {\em Automatic Control,
  IEEE Transactions on}, vol.~15, no.~1, pp.~44 -- 48, 1970.

\bibitem{graphtheory}
D.~B. West, {\em Introduction to Graph Theory}.
\newblock Prentice Hall, 2001.

\bibitem{Handbook1996}
H.~L\"{u}tkepohl, {\em Handbook of matrices}.
\newblock Wiley, 1996.

\bibitem{O'Reilly1981363}
J.~O'Reilly, ``The discrete linear time invariant time-optimal control
  problem--an overview,'' {\em Automatica}, vol.~17, no.~2, pp.~363 -- 370,
  1981.

\bibitem{Siljak1991}
D.~D. {\v{S}}iljak, {\em Decentralized control of complex systems}.
\newblock No.~184, Academic Press, 1991.

\bibitem{Farokhi2011}
F.~Farokhi and K.~H. Johansson, ``Dynamic control design based on limited model
  information,'' in {\em Communication, Control, and Computing, 49th Annual
  Allerton Conference on}, pp.~1576 -- 1583, 2011.

\bibitem{Komaroff1994}
N.~Komaroff, ``Iterative matrix bounds and computational solutions to the
  discrete algebraic \textsc{R}iccati equation,'' {\em Automatic Control, IEEE
  Transactions on}, vol.~39, no.~8, pp.~1676 -- 1678, 1994.

\end{thebibliography}

\appendix

\section{Proof of Lemma~\ref{lem:0}} \label{prooflem:0}
For any plant $P=(A,B,x_0) \in \mathcal{P}$, the optimal controller $K^*(P)$ exists (because the plant is controllable since $B$ is invertible by assumption) and can be computed using the unique positive definite solution to the discrete algebraic Riccati equation
\begin{equation} \label{eqn:2}
X=A^TXA-A^TXB(I+B^TXB)^{-1}B^TXA+I.
\end{equation}
The corresponding cost is $J_P(K^*(A,B))=x_0^T(X-I)x_0$. Inserting the product $BB^{-1}$ before every matrix $A$ and $B^{-T}B^T$ after every matrix $A^T$ in~(\ref{eqn:2}) results in
\begin{equation} \label{eqn:3}
\begin{split}
X-&I=A^TB^{-T}B^TXBB^{-1}A\\&-A^TB^{-T}B^TXB(I+B^TXB)^{-1}B^TXBB^{-1}A.
\end{split}
\end{equation}
Naming $B^TXB$ as $Y$ simplifies~(\ref{eqn:3}) into
\begin{equation} \label{eqn:4}
X-I=A^TB^{-T}[Y-Y(I+Y)^{-1}Y]B^{-1}A.
\end{equation}
Note that $Y$ is a positive definite matrix because $X$ is positive definite and $B$ is full rank. Let us denote the right-hand side of~(\ref{eqn:4}) by $A^TB^{-T}g(Y)B^{-1}A$. Then we can make the following two claims regarding the rational function $g(\cdot)$.
\par \textit{Claim 1:} The function $y \mapsto g(y)=y/(1+y)$ is a monotonically increasing over $\mathbb{R}^+.$
\par \textit{Claim 2:} Let $Y \in S_{++}^n$ and $D$, $T$ be diagonal and unitary matrices, respectively, such that $Y=T^TDT$. Then $g(Y)=T^T \diag(g(d_{ii}))T$, where $d_{ii}$ are the diagonal elements of $D$ (and the eigenvalues of $Y$).
\par Claim 1 is proved by computing the derivative of $g$ over $\mathbb{R}^+$, while Claim 2 follows from the fact that all matrices involved in the computation of $g(Y)$ can be diagonalized in the same basis yielding
\begin{equation*}
\begin{split}
g(Y)&=Y-Y(I+Y)^{-1}Y\\&=T^TDT-T^TDT(I+T^TDT)^{-1}T^TDT\\&=T^T(D-D(I+D)^{-1}D)T\\&=T^Tg(D)T.
\end{split}
\end{equation*}
Using these two claims, we find that, for all $Y$ with eigenvalues denoted by $\lambda_1(Y), \dots, \lambda_n(Y)$
\begin{equation} \label{eqn:5}
\begin{split}
X-I&=A^TB^{-T}g(Y)B^{-1}A\\&= A^TB^{-T}T^T \diag(g(\lambda_{i}(Y)))TB^{-1}A\\&\geq (g(\underline{\lambda}(Y))) A^TB^{-T}B^{-1}A,
\end{split}
\end{equation}
where $\underline{\lambda}(Y)$ is a positive number because matrix $Y$ is a positive definite matrix. Now, according to~\cite{Komaroff1994},
\begin{equation} \label{eqn:5.5}
\underline{\lambda}(X) \geq \underline{\lambda}(A^T(I+BB^T)^{-1}A+I)\geq \frac{\underline{\sigma}^2(A)}{1+\bar{\sigma}^2(B)}+1.
\end{equation}
Using~(\ref{eqn:5.5}) in inequality $\underline{\lambda} (Y) \geq \underline{\sigma}^2(B) \underline{\lambda}(X)$ gives
\begin{equation} \label{eqn:5.6}
\underline{\lambda}(Y) \geq \frac{\underline{\sigma}^2(B)\underline{\sigma}^2(A)}{1+\bar{\sigma}^2(B)}+\underline{\sigma}^2(B),
\end{equation}
and, because of the claim~1 and the inequality in~(\ref{eqn:5.6}), we will have
\begin{equation} \label{eqn:6}
\begin{split}
g(\underline{\lambda}(Y)) &\geq \frac{\underline{\sigma}^2(B)[\underline{\sigma}^2(A)+\bar{\sigma}^2(B)+1]}
{1+\bar{\sigma}^2(B)+\underline{\sigma}^2(B)[\underline{\sigma}^2(A)+\bar{\sigma}^2(B)+1]} \\ &\geq \frac{\underline{\sigma}^2(B)}{\underline{\sigma}^2(B)+1}.
\end{split}
\end{equation}
Combining~(\ref{eqn:5}) and~(\ref{eqn:6}) results in
\begin{equation*}
X-I \geq \frac{\underline{\sigma}^2(B)}{\underline{\sigma}^2(B)+1} A^TB^{-T}B^{-1}A,
\end{equation*}
and, therefore
\begin{equation*}
\begin{split}
J_P(K^*(A,B))&=x_0^T(X-I)x_0 \\ &\geq \left( \frac{\underline{\sigma}^2(B)}{\underline{\sigma}^2(B)+1} \right) x_0^T(A^TB^{-T}B^{-1}A)x_0 \\&= \left( \frac{\underline{\sigma}^2(B)}{\underline{\sigma}^2(B)+1} \right) J_P(\Gamma^\Delta(A,B)).
\end{split}
\end{equation*}

\section{Proof of Lemma~\ref{lem:1}} \label{prooflem:1}

Let $\Gamma \in \mathcal{C}$ and assume that the implication does not hold, i.e., that there exists a matrix $A$ and indices $i$, $j$ with $\ell_0 \in \mathcal{I}_i$ such that $a_{\ell j} = 0$ for all $\ell \in \mathcal{I}_i$ but $\gamma_{\ell_0j}(A,B) \neq 0$. Consider matrix $\bar{A}$ such that $\bar{A}_\ell=A_\ell$ for all $\ell\in \mathcal{I}_i$ and $\bar{A}_z=0$ for all $z\notin \mathcal{I}_i$. Based on the definition of limited-model-information control design methods, we know that $\Gamma_\ell(\bar{A},B)=\Gamma_\ell(A,B)$ for all $\ell\in \mathcal{I}_i$ and $\Gamma_z(\bar{A},B)=0$ for all $z\notin \mathcal{I}_i$ (because $\Gamma_z(A,B) = \Gamma_z(0,B)$ for all $z\notin \mathcal{I}_i$ and, as shown in~\cite{langbort10}, it is necessary that $\Gamma(0,B)=0$ for a finite competitive ratio). For $x=e_j$, we have
\begin{equation*}
\begin{split}
J_{(\bar{A},B,e_j)}(\Gamma(\bar{A},B)) \geq \sum_{\ell\in \mathcal{I}_i} \gamma_{\ell j}(\bar{A},B)^2 &= \sum_{\ell\in \mathcal{I}_i} \gamma_{\ell j}(A,B)^2 \\&\geq \gamma_{\ell_0j}(A,B)^2 > 0.
\end{split}
\end{equation*}
Using~(\ref{inequalityusingdeadbeat}), we get
\begin{equation*}
\begin{split}
r_{\p}(\Gamma) \geq \frac{J_{(\bar{A},B,e_j)}(\Gamma(\bar{A},B))}{J_{(\bar{A},B,e_j)}(\Gamma^\Delta(\bar{A},B))} = \infty,
\end{split}
\end{equation*}
since $J_{(\bar{A},B,e_j)}(\Gamma^\Delta(\bar{A},B))=0$. This proves the claim by contrapositive.

\section{Proof of Lemma~\ref{lem:2}} \label{prooflem:2}
Clearly, it is enough to prove inequality~(\ref{eq_rat}) for control design methods with a finite competitive ratio.

We proceed in three steps. First, using Lemma~\ref{lem:1}, we characterize the design strategies leading to a finite competitive ratio. Then, we argue that the controllers produced by such strategies must be stabilizing for all plants, and use the fact that every closed-loop characteristic polynomial is Schur to construct a sequence of real numbers with specific properties for each control design strategy. We then use this sequence to construct a sequence of plants allowing us to lower bound the competitive ratio of each control design strategy.

Let $G_{\mathcal{P}}$ have a loop and $\Gamma \in \mathcal{C}$ have finite competitive ratio. Without loss of generality, let us assume that the nodes of graph $G_{\mathcal{P}}$ are numbered such that it admits the following loop of length $\ell$: $1 \rightarrow 2 \rightarrow \cdots \rightarrow \ell \rightarrow 1$. Let us choose indices $i_1 \in \mathcal{I}_1$, $i_2 \in \mathcal{I}_2$, $\ldots$, $i_\ell \in \mathcal{I}_\ell$ and consider the one-parameter family of matrices $\{A(r)\}$ defined by $a_{i_2i_1}(r)=r$, $a_{i_3i_2}(r)=r$, $\ldots$, $a_{i_\ell i_{\ell-1}}(r)=r$, $a_{i_1i_\ell}(r)=r$, and all other entries equal to zero, for all $r$. Let $B=\epsilon I$. Because of Lemma~\ref{lem:1}, the controller gain entries $\gamma_{j_2i_1}(A(r),B)$ for all $j_2 \in \mathcal{I}_2$, $\gamma_{j_3i_2}(A(r),B)$ for all $j_3 \in \mathcal{I}_3$, $\ldots$, $\gamma_{j_\ell i_{\ell-1}}(A(r),B)$ for all $j_\ell \in \mathcal{I}_\ell$, $\gamma_{j_1i_\ell}(A(r),B)$ for all $j_1 \in \mathcal{I}_1$ can be non-zero, but all other entries of the controller gain $\Gamma(A(r),B)$ are zero for all $r$. As a result, the characteristic polynomial of matrix $A(r)+B\Gamma(A(r),B)$ can be computed as:
\begin{equation*}
\begin{split}
\lambda^{n-\ell}[\lambda^\ell-&(-1)^\ell(r+\epsilon\gamma_{i_2i_1}(A(r),B))(r+\epsilon\gamma_{i_3i_2}(A(r),B))\\ \times \cdots \times &(r+\epsilon\gamma_{i_\ell i_{\ell-1}}(A(r),B))(r+\epsilon\gamma_{i_1i_\ell}(A(r),B))].
\end{split}
\end{equation*}
Now, note that because $\Gamma$ has a bounded competitive ratio against $\mathcal{P}$ by assumption, this polynomial should be stable for all $r$. (Indeed, $\Gamma$ can have a finite competitive ratio only if $A+B\Gamma(A,B)$ is stable for all matrices $A$, otherwise it would yield an infinite cost for some plants while the corresponding optimal cost remains bounded since the pair $(A,B)$ is controllable for all plant in $\mathcal{P}$). As a result, we must have
\begin{equation}
\label{first_eq}
\begin{split}
&\left| (r+\epsilon\gamma_{i_2i_1}(A(r),B))\cdots(r+\epsilon\gamma_{i_1i_\ell}(A(r),B)) \right| \\& \hspace{0.10in} =|r+\epsilon\gamma_{i_2i_1}(A(r),B)| \cdots |r+\epsilon\gamma_{i_1i_\ell}(A(r),B)|<1
\end{split}
\end{equation}
for all $r$. Let $\{r_z\}_{z=1}^\infty$ be a sequence of real numbers with the property that $r_z$ goes to infinity as $z$ goes to infinity. From~(\ref{first_eq}), we know that there exists an index $\bar{m}$ such that
\begin{equation}
\label{big_prop}
\forall N,\; \exists z >N \mbox{ such that } |r_z+\epsilon \gamma_{i_{\bar{m} \oplus 1}i_{\bar{m}}}(A(r_z),B)|<1,
\end{equation}
where ``$\oplus$'' designated addition modulo $\ell$; i.e., $i\oplus j=(i+j)-\lfloor(i+j)/\ell\rfloor \ell$ where $\lfloor x \rfloor=\max\{y\in \mathbb{Z}|y\leq x \}$ for all $x\in \mathbb{R}$. Indeed, if this is not the case, it is true that
\begin{equation*}
\begin{split}
\forall m, \; \exists N_m \mbox{ such that } |r_z+\epsilon \gamma_{i_{m \oplus 1}i_{m}}(A(r_z),B)&| \geq 1, \\ &\forall z>N_m.
\end{split}
\end{equation*}
Then, for all $z > \max_{m} N_m$ and all $m$, $$|r_z+\epsilon \gamma_{i_{m \oplus 1}i_{m}}(A(r_z),B)| \geq 1$$ which contradicts~(\ref{first_eq}). Without loss of generality (since this just amounts to renumbering the nodes in the plant graph), we assume that $\bar{m} = 1$. Using~(\ref{big_prop}), we can then construct a subsequence $\{r_{\phi(z)}\}$ of $\{r_z\}$ with the property that
$$
|r_{\phi(z)}+\epsilon \gamma_{i_{2}i_1}(A(r_{\phi(z)}),B)|< 1 \; \mbox{ for all } z.
$$
Now introduce the sequence of matrices $\{\bar{A}(z)\}_{z=1}^{\infty}$ defined by $\bar{A}_{i_2i_1}(z)=r_{\phi(z)}$ for all $z$ and every other row equal to zero. For large enough $z$ (and hence, large enough $r_{\phi(z)}$), we get
\begin{equation*}
\begin{split}
J_{(\bar{A}(z),B,e_{i_1})}(\Gamma(\bar{A}(z),B)) &\geq \gamma_{i_2i_1}(\bar{A}(z),B)^2 \\ &= \gamma_{i_2i_1}(A(r_{\phi(z)}),B)^2 \\ &\geq \frac{(|r_{\phi(z)}|-1)^2}{\epsilon^2},
\end{split}
\end{equation*}
and thus
\begin{equation*}
\frac{J_{(\bar{A}(z),B,e_{i_1})}(\Gamma(\bar{A}(z),B))}{J_{(\bar{A}(z),B,e_{i_1})}(K^*(\bar{A}(z),B,e_{i_1}))}\geq \frac{(|r_{\phi(z)}|-1)^2/\epsilon^2}{r_{\phi(z)}^2/(1+\epsilon^2)}.
\end{equation*}
This, in particular, implies that
\begin{equation*}
r_{\p}(\Gamma) \geq \lim_{z \rightarrow \infty}\frac{J_{(\bar{A}(z),B,e_{i_1})}(\Gamma(\bar{A}(z),B))}{J_{(\bar{A}(z),B,e_{i_1})}(K^*(\bar{A}(z),B,e_{i_1}))}\geq 1+1/\epsilon^2.
\end{equation*}
 Note that $\bar{A}(z)$ is a nilpotent matrix for all $z$, and thus
$$
J_{(\bar{A}(z),B,e_{i_1})}(K^*(\bar{A}(z),B,e_{i_1}))=J_{(\bar{A}(z),B,e_{i_1})}(K^*_C(\bar{A}(z),B))
$$
similar to the proof of Theorem~\ref{tho:1}, and therefore
$$
J_{(\bar{A}(z),B,e_{i_1})}(K^*_C(\bar{A}(z),B))=r_{\phi(z)}^2/(1+\epsilon^2)
$$
using the unique positive-definite solution of discrete algebraic Riccati equation in~(\ref{eqn:Riccati}).

\section{Proof of Lemma~\ref{lem:undomination}} \label{prooflem:undomination}

We prove that if there is no sink in the plant graph (i.e., according to~\cite{graphtheory}, if $\forall j \exists k$, $k\neq j$, such that $(s_{\mathcal{P}})_{kj}\neq 0$) then the deadbeat control design method is undominated. For proving this claim, we are going to prove that for any control design $\Gamma\in \comm \backslash \{\ \Gamma^\Delta \}$, there exists a plant $P=(A,B,x_0)\in \p$ such that $J_P(\Gamma(A,B))>J_P(\Gamma^\Delta(A,B))=x_0^T[A^TB^{-T}B^{-1}A]x_0$. We will proceed in several steps, which require us to partition the set of limited model information control design methods $\comm$ as follows
\begin{equation*}
\comm=\mathcal{L}^c \cup \mathcal{W}_1 \cup \mathcal{W}_2 \cup \{\Gamma^\Delta\},
\end{equation*}
where
\begin{equation*}
\begin{split}
&\mathcal{L}:= \{\Gamma \in \comm | \exists \Lambda_j:\Re^{n_j \times n}\times \Re^{n_j \times n_j} \rightarrow \Re^{n_j \times n_j}, \\ [\Gamma &(A,B)]_j = \Lambda_j([A]_j,B_{jj}) [A]_j,\hspace{0.05in} \textrm{for all} \hspace{0.05in} j=1,\cdots,q \},
\end{split}
\end{equation*}
\begin{equation*}
\begin{split}
\mathcal{W}_1&:= \{\Gamma\in \mathcal{L} | \exists j, i\neq j \hspace{0.05in} \textrm{and} \hspace{0.05in} A_{ij} \in \Re^{n_i\times n_j} \hspace{0.05in} \textrm{nonzero s.t.} \hspace{0.05in}\\
&I+B_{ii}\Lambda_i(\left[ 0 \; \cdots \; 0 \; A_{ij} \; 0 \; \cdots \; 0 \right],B_{ii})\neq 0\},
\end{split}
\end{equation*}
and
\begin{equation*}
\begin{split}
\mathcal{W}_2&:= \{\Gamma\in \mathcal{L}\setminus \mathcal{W}_1 | \exists i\in \{1,\cdots,q \}, [A]_i\in \Re^{n_i \times n}, \hspace{0.05in} \textrm{with} \\ & \textrm{appropriate structure s.t.} \hspace{0.05in} I+B_{ii}\Lambda_i([A]_i,B_{ii})\neq 0 \}.
\end{split}
\end{equation*}
 In words, $\mathcal{L}$ is the set of all control design methods for which sub-controller $K_j$ can be written as a linear combination of vectors in $\{A_i , i \in \mathcal{I}_j \}$ for all $j$. Sets $\mathcal{W}_1$ and $\mathcal{W}_2$ are subsets of $\mathcal{L}$ which put further constraints on map $\Gamma$. Using different lower bounds on closed-loop performance in each case, we show that $\Gamma^{\Delta}$ is undominated by control strategies in each of $\mathcal{L}^c$, $\mathcal{W}_1$, and $\mathcal{W}_2$.

First, we prove that the deadbeat control design method is undominated by control design strategies in $\mathcal{L}^c$. Let $\Gamma \in \mathcal{L}^c$ and let $j$ be such that $\exists j_1\in \mathcal{I}_j$ which $\Gamma_{j_1}(\bar{A},B)^T$ cannot be written as a linear combination of vectors in the set $\{\bar{A}_{i}^T, \forall i\in \mathcal{I}_j \}$ for some matrix $\bar{A}$ and matrix $B$. Let $a_{i}^T=\bar{A}_{i}$ for all $i\in \mathcal{I}_j$ and consider matrix $A$ such that the row $A_{i}=a_{i}^T$ for all $i\in \mathcal{I}_j$ and $A_i=0$ for all $i \in \mathcal{I}_j^c$. If $\Gamma(0,B)\neq 0$, then $\Gamma$ cannot dominate $\Gamma^\Delta$ (since $\Gamma^\Delta(0,B)=0$ for all $x_0$) and, thus, there is no loss of generality in assuming that $\Gamma(0,B)=0$ for all $x_0$, and, in turn that $\Gamma_i(A,B)=0$ for all $i\in \mathcal{I}_j^c$. Let us also denote $\Gamma(A,B)$ by $K$ and $\Gamma_i(A,B)=\Gamma_i(\bar{A},B)$ by $K_i^T$ for all $i\in \mathcal{I}_j$. For all $x_0$,
\begin{equation*}
\begin{split}
J_{(A,B,x_0)}&(\Gamma(A,B)) \geq x_0^T[K^TK+(A+BK)^T(A+BK)]x_0,
\end{split}
\end{equation*}
and
\begin{equation} \label{eqn:7}
\begin{split}
J_{(A,B,x_0)}(\Gamma(A,B)&)-J_{(A,B,x_0)}(\Gamma^\Delta(A,B)) \\ \geq x_0^T[A^T&(I-B^{-T}B^{-1})A+A^TBK\\+K^T&B^TA+K^T(I+B^TB)K]x_0.
\end{split}
\end{equation}
We know that $\nullf (A) = \spanf \{A_{i}^T, \empty \forall i\in \mathcal{I}_j \}^\bot \neq \{0 \}$, because $n_j < n$. On the other hand, we know that there exists an $j_1 \in \mathcal{I}_j$ such that $K_{j_1}\notin \spanf \{A_{i}^T, \empty \forall i\in \mathcal{I}_j \}$ which shows that
\begin{equation*}
\begin{split}
\spanf \{A_{i}^T, \empty \forall i\in \mathcal{I}_j \} \varsubsetneq \spanf \{A_{i}^T, \empty &\forall i\in \mathcal{I}_j \} +\spanf \{K_{i}^T, \empty \forall i\in \mathcal{I}_j \},
\end{split}
\end{equation*}
Thus, we can choose an initial condition $x_0 \in \nullf (A)$ such that $Kx_0 \neq 0$. Using this $x_0$ in~(\ref{eqn:7}) results in
\begin{equation}
\begin{split}
J_{(A,B,x_0)}(\Gamma(A,B))&-J_{(A,B,x_0)}(\Gamma^\Delta(A,B)) \geq x_0^T[K^T(I+B^TB)K]x_0 > 0.
\end{split}
\end{equation}
Therefore, the control design strategies in $\mathcal{L}^c$ cannot dominate the deadbeat control design strategy $\Gamma^\Delta$.
\par Second, we prove that the deadbeat control design strategy is undominated by control design methods in $\mathcal{W}_1$. Let $\Gamma\in \mathcal{W}_1$ and let $j$ be such that $(I+B_{ii}\Lambda_i(\left[0 \; \cdots \; 0 \; \bar{A}_{ij} \; 0 \; \cdots \; 0\right],B_{ii}))\neq 0$ for some $i\neq j$. It means that there exists at least $i_1\in \mathcal{I}_i$ and $j_1\in \mathcal{I}_j$ such that $\bar{a}_{i_1j_1}\neq 0$ and $\bar{a}_{i_1j_1}+b_{i_1i_1}\gamma_{i_1j_1}(\bar{A},B)\neq 0$. Using the structure matrix, we know that there exits a $\ell\neq i$ such that $(s_{\mathcal{P}})_{\ell i}\neq 0$. Choose an index $\ell_1\in \mathcal{I}_\ell$. Consider the matrix $A$ defined by $[A]_i=[\bar{A}]_i$, $a_{\ell_1i_1}=r$ and all other entries equal to zero. Then, $[\Gamma(A,B)]_i=\Lambda_i([A]_i,B_{ii})[A]_i$, $[\Gamma(A,B)]_\ell=\Lambda_\ell([A]_\ell,B_{\ell\ell})[A]_\ell$ (because $\Gamma \in \mathcal{L}$), and $[\Gamma(A,B)]_z=0$ for all $z\neq i,\ell$. Denote $\Gamma(A,B)$ by $K$. We have
\begin{equation*}
\begin{split}
J_{(A,B,x_0)}(\Gamma(A,B))\geq x_0^T[(A+BK)^T&K^TK(A+BK)\\+(&(A+BK)^2)^T(A+BK)^2]x_0.
\end{split}
\end{equation*}
Using $x_0=e_{j_1}$ results in
\begin{equation} \label{eqn:8}
\begin{split}
J_{(A,B,e_{j_1})}&(\Gamma(A,B))-J_{(A,B,e_{j_1})}(\Gamma^\Delta(A,B))\geq \\ &[k_{\ell_1i_1}^2+(r+b_{\ell_1\ell_1}k_{\ell_1i_1})^2](a_{i_1j_1}+b_{i_1i_1}k_{i_1j_1})^2-\sum_{z \in \mathcal{I}_i}\frac{a_{zj_1}^2}{b_{zz}^2}.
\end{split}
\end{equation}
Note that, irrespective of the choice of the controller gain $k_{\ell_1i_1}$,
\begin{equation*}
k_{\ell_1i_1}^2+(r+b_{\ell_1\ell_1}k_{\ell_1i_1})^2 \geq r^2/(1+b_{\ell_1\ell_1}^2),
\end{equation*}
and as a result,
\begin{equation*}
\lim_{r\rightarrow +\infty} [k_{\ell_1i_1}^2+(r+b_{\ell_1\ell_1}k_{\ell_1i_1})^2](a_{i_1j_1}+b_{i_1i_1}k_{i_1j_1})^2=\infty,
\end{equation*}
because $a_{i_1j_1}+b_{i_1i_1}k_{i_1j_1}\neq 0$. Hence, we can always construct $A$ with appropriate choice of index $\ell$ and a scalar $r$ large enough to make the right hand side of the expression~(\ref{eqn:8}) positive. As a result, $\Gamma\in \mathcal{W}_1$ cannot dominate $\Gamma^\Delta$.
\par Third, we prove that the deadbeat control design strategy is undominated by control design methods in $\mathcal{W}_2$. Let $\Gamma\in \mathcal{W}_2$ and index $i$ and vector $[\bar{A}]_i$ be such that $I+\Lambda_i([\bar{A}]_i,B_{ii})\neq 0$. Thus we know that there exists at least $i_1\in \mathcal{I}_i$ such that $\bar{A}_{i_1} \neq 0$ and $\bar{A}_{i_1}+b_{i_1i_1}\Gamma_{i_1}(\bar{A},B) \neq 0$. Based on the structure matrix we know that there exits $\ell\neq i$ such that $(s_{\mathcal{P}})_{\ell i}\neq 0$. Choose an index $\ell_1\in \mathcal{I}_\ell$. Consider the matrix $A$ defined by $[A]_i=[\bar{A}]_i$ and $a_{\ell_1i_1}=r$ and all other entries of $A$ equal to zero. Then $[A]_i+B_{ii}[\Gamma(A,B)]_i=(I+B_{ii}\Lambda_i([A]_i,B_{ii}))[A]_i$ and $[A]_j+B_{jj}[\Gamma(A,B)]_j=0$ for all $j\neq i$ (and, in particular, $j=\ell$ since $\Gamma$ does not belong to $\mathcal{W}_1$). Again, $K$ will stand for $\Gamma(A,B)$. We have
\begin{equation*}
\begin{split}
K^TK+(A+BK)^T&K^TK(A+BK)-A^TB^{-T}B^{-1}A\\ \geq (A_{i_1}+b_{i_1i_1}&\Gamma_{i_1}(A,B))^T(A_{i_1}+b_{i_1i_1}\Gamma_{i_1}(A,B))\times \\& r^2/b_{\ell_1\ell_1}^2 -\sum_{z\in \mathcal{I}_i} A_z^TA_z/b_{zz}^2,
\end{split}
\end{equation*}
and hence, since $A_{i_1}+b_{i_1i_1}\Gamma_{i_1}(A,B) \neq 0$, we can choose $r$ large enough to ensure that this matrix has a strictly positive eigenvalue. Thus, the control design strategy $\Gamma\in \mathcal{W}_2$ cannot dominate $\Gamma^\Delta$.

\section{Proof of Lemma~\ref{lem:3}} \label{prooflem:3}

The proof is by contrapositive. Let $\Gamma$ be communication-less and assume that there exist matrices $A$ and $B$ and indices $i_1\in \mathcal{I}_i$ and $j_1\in \mathcal{I}_j$ such that $a_{i_1j_1}+b_{i_1i_1}\gamma_{i_1j_1}(A,B)\neq 0$. Choose an index $k_1\in \mathcal{I}_k$. Consider the one-parameter family of matrices $\bar{A}(r)$ defined by $[\bar{A}(r)]_i=[A]_i$, $\bar{a}_{k_1i_1}=r$, and all other entries of $\bar{A}(r)$ being equal to zero for all $r$. We know that $[\Gamma(\bar{A}(r),B)]_i=[\Gamma(A,B)]_i$ and $\Gamma_{\bar{k}}(\bar{A}(r),B)=\gamma_{\bar{k}i_1}(r)e_{i_1}^T$ for all $\bar{k} \in \mathcal{I}_k$ (because of Lemma~\ref{lem:1}), $[\Gamma(\bar{A}(r),B)]_z=0$ for all $z\neq i,k$. For $x_0=e_{j_1}$, we have
\begin{equation*}
\begin{split}
J_{(\bar{A}(r),B,e_{j_1})}(\Gamma(\bar{A}(r),B)) \geq (&a_{i_1j_1}+b_{i_1i_1}\gamma_{i_1j_1}(A,B))^2\\&\times[\gamma_{k_1i_1}(r)^2+(r+b_{k_1k_1}\gamma_{k_1i_1}(r))^2].
\end{split}
\end{equation*}
The minimum value of function $y \mapsto [y^2+(r+b_{k_1k_1}y)^2]$ is $r^2/(1+b_{k_1k_1}^2)$. Hence, irrespective of function $\gamma_{k_1i_1}$,
\begin{equation*}
\begin{split}
J_{(\bar{A}(r),B,e_{j_1})}(\Gamma(\bar{A}(r),B)) \geq (a_{i_1j_1}+&b_{i_1i_1}\gamma_{i_1j_1}(A,B))^2 r^2/(1+b_{k_1k_1}^2).
\end{split}
\end{equation*}
Note that the term $(a_{i_1j_1}+b_{i_1i_1}\gamma_{i_1j_1}(A,B))^2$ is independent from $r$ because $\Gamma$ is communication-less. In addition,
\begin{equation*}
J_{(\bar{A}(r),B,e_{j_1})}(\Gamma^\Delta(\bar{A}(r),B))=\sum_{z\in \mathcal{I}_i} \frac{\bar{a}_{zj_1}^2}{b_{zz}^2}=\sum_{z\in \mathcal{I}_i} \frac{a_{zj_1}^2}{b_{zz}^2}
\end{equation*}
for all $r$ and, thus, $J_{(\bar{A}(r),B,e_{j_1})}(\Gamma^\Delta(\bar{A}(r),B))$ is also independent from $r$. Then, proceeding as in~(\ref{inequalityusingdeadbeat}), we deduce that
\begin{equation*}
r_\p(\Gamma) \geq \frac{(a_{i_1j_1}+b_{i_1i_1}\gamma_{i_1j_1}(A,B))^2}{(1+b_{k_1k_1}^2)J_{(\bar{A}(r),B,e_{j_1})}(\Gamma^\Delta(\bar{A}(r),B))} \lim_{r\rightarrow \infty} r^2.
\end{equation*}
Since $(a_{i_1j_1}+b_{i_1i_1}\gamma_{i_1j_1}(A,B)) \neq 0$ by assumption, we then deduce that $\Gamma$ has an unbounded competitive ratio, which proves the lemma by contrapositive.

\section{Proof of Theorem~\ref{tho:4}} \label{prooftho:4}

We prove that for any control design method $\Gamma\in \comm \backslash \{ \Gamma^\Theta \}$, there exists a plant $P=(A,B,x_0)\in \p$ such that $J_P(\Gamma(A,B))>J_P(\Gamma^\Theta(A,B))$. Like in the proof of Theorem 3.6, we partition the set of limited model information control design methods $\comm$ as follows
\begin{equation*}
\mathcal{C}=\mathcal{L}^c \cup \mathcal{W}_0 \cup \mathcal{W}_1 \cup \mathcal{W}_2 \cup \{\Gamma^\Theta\},
\end{equation*}
where
\begin{equation*}
\begin{split}
\mathcal{L}:= \{\Gamma \in \comm | & \exists \Lambda_i:\Re^{n_i \times n}\times \Re^{n_i \times n_i} \rightarrow \Re^{n_i \times n_i}, \\ [\Gamma &(A,B)]_i = \Lambda_i([A]_i,B_{ii}) [A]_i,\hspace{0.05in} \textrm{for all} \hspace{0.05in} i=1,\cdots,q \},
\end{split}
\end{equation*}
\begin{equation*}
\begin{split}
\mathcal{W}_0:= \{\Gamma\in& \mathcal{L}, \exists i \in \{ q-c+1,\dots,q\} \hspace{0.05in} \textrm{such that} \hspace{0.05in} \Lambda_i([A]_i,B_{ii})\neq W_i([A]_i,B_{ii}) \hspace{0.05in}\},\\
\end{split}
\end{equation*}
with $W_i$ defined as in equation~(\ref{eqn:12}),
\begin{equation*}
\begin{split}
\mathcal{W}_1:= \{\Gamma\in \mathcal{L} \setminus \mathcal{W}_0| \exists &i\in \{1,\cdots,q-c \}, \exists j\neq i \hspace{0.05in} \textrm{and} \; A_{ij} \in \Re^{n_i\times n_j} \hspace{0.05in} \textrm{nonzero} \\& \mbox{such that}
I+B_{ii}\Lambda_i(\left[ 0 \; \cdots \; 0 \; A_{ij} \; 0 \; \cdots \; 0 \right],B_{ii})\neq 0\},\\
\end{split}
\end{equation*}
and
\begin{equation*}
\begin{split}
\mathcal{W}_2&:= \{\Gamma\in \mathcal{L}\setminus \mathcal{W}_0 \cup \mathcal{W}_1 | \exists i\in \{1,\cdots,q-c \}, \;[A]_i\in \Re^{n_i \times n},\\ &\textrm{ with appropriate structure such that} \hspace{0.05in} I+B_{ii}\Lambda_i([A]_i,B_{ii})\neq 0 \}.
\end{split}
\end{equation*}
First, we prove that $\Gamma^\Theta$ is undominated by control design methods in $\mathcal{L}^c$. Let $\Gamma\in \mathcal{L}^c$ and let $i$ be such that there exists a plant with matrix $\bar{A}$ with the property that subcontroller $[\Gamma]_i([\bar{A}]_i,B_{ii})^T$ does not belong to the linear subspace spanned by the columns of $[\bar{A}]_i^T$. If $1 \leq i \leq q-c$ then, proceeding as in the proof of Theorem~\ref{tho:2}, we can find matrices $A$, $B$ and initial condition $x_0$ such that $J_{P}(\Gamma(P)) > J_{P}(\Gamma^{\Delta}(P)) = J_{P}(\Gamma^{\Theta}(P))$ for $P=(A,B,x_0)$ (with the last equality following from the structure of matrix $A$). Hence, without loss of generality, we assume that $q-c+1 \leq i \leq q$. Consider matrix $A$ defined as $[A]_i=[\bar{A}]_i$ and $[A]_j=0$ for all $j\neq i$. For this particular matrix $A$ and any $B$, $x_0$ we know from the proof of the ``only if'' part of the Theorem~\ref{tho:2} that $\Gamma^\Theta(A,B,x_0)$ is the globally optimal controller. Hence, every other control design method in $\comm$ leads to a controller with greater performance criterion than $\Gamma^\Theta$ for this particular type of plants. Therefore, the control design $\Gamma^\Theta$ is undominated by control design methods in $\mathcal{L}^c$.

The same reasoning shows that $\Gamma^{\Theta}$ is also undominated by control design methods in $\mathcal{W}_0$.

We now prove that $\Gamma^\Theta$ is undominated by control design strategies in $\mathcal{W}_1$. Let $\Gamma\in \mathcal{W}_1$ and let $1 \leq i\leq q-c$ be such that $(I + B_{ii} \Lambda_i ( [\bar{A}]_i , B_{ii} ) ) \neq 0$ where $[\bar{A}]_i=\left[0 \; \cdots \; 0 \; \bar{A}_{ij} \; 0 \; \cdots \; 0\right]$ for some $j\neq i$. This means that there exists at least one $i_1\in \mathcal{I}_i$ and $j_1\in \mathcal{I}_j$ such that $\bar{a}_{i_1j_1}\neq 0$ and $\bar{a}_{i_1j_1}+b_{i_1i_1}\gamma_{i_1j_1}(A,B)\neq 0$. Because subsystem $i$ is not a sink (since $1\leq i\leq q-c$), we know that there exists a $z\neq i$ such that $(s_{\mathcal{P}})_{zi}\neq 0$. If $1 \leq z\leq q-c$ we can again proceed as in the proof of Theorem~\ref{tho:2} to construct a plant $P$ for which $J_{P}(\Gamma(P)) > J_{P}(\Gamma^{\Theta}(P))$. Thus, without loss of generality, we assume that $q-c+1 \leq z\leq q$. Choose an index $z_1\in \mathcal{I}_z$ and consider the matrix $A$ defined by $[A]_i=[\bar{A}]_i$, $a_{z_1i_1}=r$ and all other entries equal to zero. Then, $[\Gamma(A,B)]_i=\Lambda_i([A]_i,B_{ii})[A]_i$, $[\Gamma(A,B)]_z=-b_{z_1z_1}/(1+b_{z_1z_1}^2)[A]_z$ (because $\Gamma \notin \mathcal{W}_0 \cup \mathcal{L}^c$), and $[\Gamma(A,B)]_t=0$ for all $t\neq i,z$. Denoting $\Gamma(A,B)$ by $K$, we see that
\begin{equation*}
\begin{split}
J_{(A,B,x_0)}(\Gamma(A,B))\geq x_0^T[(A+BK)^T&K^TK(A+BK)\\+(&(A+BK)^2)^T(A+BK)^2]x_0
\end{split}
\end{equation*}
for all $B \in \mathcal{B}(\epsilon)$ and $x_0$. Taking $x_0=e_{j_1}$ then results in
\begin{equation}
\begin{split}
J_{(A,B,e_{j_1})}&(\Gamma(A,B))-J_{(A,B,e_{j_1})}(\Gamma^\Theta(A,B))\geq \\ &[k_{z_1i_1}^2+(r+b_{z_1z_1}k_{z_1i_1})^2](a_{i_1j_1}+b_{i_1i_1}k_{i_1j_1})^2-\sum_{t \in \mathcal{I}_i}\frac{a_{tj_1}^2}{b_{tt}^2}.
\end{split}
\end{equation}
Note that, irrespective of the choice of the controller gain $k_{z_1i_1}$,
\begin{equation*}
k_{z_1i_1}^2+(r+b_{z_1z_1}k_{z_1i_1})^2 \geq r^2/(1+b_{z_1z_1}^2),
\end{equation*}
and as a result,
\begin{equation*}
\lim_{r\rightarrow +\infty} [k_{z_1i_1}^2+(r+b_{z_1z_1}k_{z_1i_1})^2](a_{i_1j_1}+b_{i_1i_1}k_{i_1j_1})^2=+\infty,
\end{equation*}
because $a_{i_1j_1}+b_{i_1i_1}k_{i_1j_1}\neq 0$. Hence, we can always construct $A$ with appropriate choice of index $z$ and a scalar $r$ large enough to make the cost difference positive. As a result, $\Gamma$ cannot dominate $\Gamma^\Theta$.
\par Finally, we prove that $\Gamma^\Theta$ is undominated by control design methods in $\mathcal{W}_2$. Let $\Gamma\in \mathcal{W}_2$ and index $1\leq i\leq q-c$ and model sub-matrices $[\bar{A}]_i$ and $B_{ii}$ such that $I+\Lambda_i([\bar{A}]_i,B_{ii})\neq 0$. Therefore, we know that there exists at least one index $i_1\in \mathcal{I}_i$ such that $\bar{A}_{i_1} \neq 0$ and $\bar{A}_{i_1}+b_{i_1i_1}\Gamma_{i_1}(\bar{A},B) \neq 0$. Based on the fact that node $i$ is not a sink, we know that there exists $z\neq i$ such that $(s_{\mathcal{P}})_{zi}\neq 0$. For the same reasons as before we again restrict ourselves, without loss of generality, to the case where $q-c+1 \leq z\leq q$. Consider the matrix $A$ defined by $[A]_i=[\bar{A}]_i$ and $a_{z_1i_1}=r$ and all other entries of $A$ equal to zero. Then, $[A]_i+[\Gamma(A,B)]_i=(I+\Lambda_i([A]_i,B_{ii}))[A]_i$ and $[\Gamma(A,B)]_z=-b_{z_1z_1}/(1+b_{z_1z_1}^2)[A]_z$ (because $\Gamma \notin \mathcal{W}_0 \cup \mathcal{L}^c$). Again, $K$ will stand for $\Gamma(A,B)$. Then, for all $B \in \mathcal{B}(\epsilon)$ and $x_0$
\begin{equation*}
\begin{split}
J_{(A,B,x_0)}(\Gamma(A,B&))-J_{(A,B,x_0)}(\Gamma^\Theta(A,B)) \\ \geq x_0^T(A_{i_1}+b_{i_1i_1}&\Gamma_{i_1}(A,B))^T(A_{i_1}+b_{i_1i_1}\Gamma_{i_1}(A,B))x_0\times \\& \hspace{-.14in} r^2b_{z_1z_1}^2/(1+b_{z_1z_1}^2)^2 -\sum_{t\in \mathcal{I}_i} x_0^TA_t^TA_tx_0/b_{tt}^2,
\end{split}
\end{equation*}
and hence, since $A_{i_1}+b_{i_1i_1}\Gamma_{i_1}(A,B) \neq 0$, we can choose $r$ large enough to ensure that this difference is strictly positive for some $x_0\in \mathbb{R}^n$ since the inner matrix will have a strictly
positive eigenvalue for large values of $r$. Thus, the control design strategy $\Gamma\in \mathcal{W}_2$ cannot dominate the control design $\Gamma^\Theta$.

\end{document}